\documentclass[10pt]{amsart}
\usepackage[pagebackref]{hyperref}

\font\teneur=eurm10

\font\sans=cmss10

\renewcommand{\d}{\mathrm{d}}

\renewcommand{\Re}{\operatorname{Re}}
\renewcommand{\Im}{\operatorname{Im}}

\newcommand{\ts}{\textstyle }
\newcommand{\cA}{{\mathcal A}}

\newcommand{\cE}{{\mathcal E}}

\newcommand{\cG}{{\mathcal G}}
\newcommand{\cH}{{\mathcal H}}
\newcommand{\cI}{{\mathcal I}}
\newcommand{\cJ}{{\mathcal J}}
\newcommand{\cK}{{\mathcal K}}
\newcommand{\cW}{{\mathcal W}}
\newcommand{\cZ}{{\mathcal Z}}
\newcommand{\bbR}{{\mathbb R}}
\newcommand{\bbB}{{\mathbb B}}
\newcommand{\bbC}{{\mathbb C}}

\newcommand{\bbP}{{\mathbb P}}

\newcommand{\bbT}{{\mathbb T}}
\newcommand{\bbZ}{{\mathbb Z}}

\newcommand{\E}{{\mathrm e}}
\newcommand{\iC}{{\mathrm i}}
\newcommand{\GL}{\operatorname{GL}}
\newcommand{\Gtwo}{\operatorname{G}_2}
\newcommand{\SL}{\operatorname{SL}}
\newcommand{\SO}{\operatorname{SO}}
\newcommand{\SU}{\operatorname{SU}}

\newcommand{\U}{\operatorname{U}}
\newcommand{\Spin}{\operatorname{Spin}}
\newcommand{\Symp}{\operatorname{Sp}}
\newcommand{\Diff}{\operatorname{Diff}}

\newcommand{\diag}{\operatorname{diag}}

\newcommand{\id}{\operatorname{id}}
\newcommand{\oppar}{{\upshape(}}
\newcommand{\clpar}{{\upshape)}}

\DeclareMathOperator{\tr}{tr}

\DeclareMathOperator{\End}{End}
\DeclareMathOperator{\Aut}{Aut}
\DeclareMathOperator{\Hom}{Hom}

\newcommand{\phm}{\phantom{-}}
\newcommand{\p}{\partial}
\newcommand{\pbar}{{\bar\partial}}
\newcommand{\bk}{{\bar k}}
\newcommand{\bj}{{\bar\jmath}}
\newcommand{\bi}{{\bar\imath}}

\newcommand{\Fs}{\mathsf{F}}

\newcommand{\Js}{\mathsf{J}}

\newcommand{\us}{\mathsf{u}}
\newcommand{\Us}{\mathsf{U}}

\newcommand{\CU}{\mathsf{CU}}
\newcommand{\QU}{\mathsf{QU}}
\newcommand{\QJ}{\mathsf{QJ}}

\newcommand{\Gams}{\hbox{\sans\char0}}
\newcommand{\psib}{\hbox{\teneur\char'040}}
\newcommand{\Phib}{\hbox{\teneur\char'010}}
\newcommand{\omegab}{\hbox{\teneur\char'041}}

\newcommand{\eugl}{\operatorname{\mathfrak{gl}}}

\newcommand{\eusu}{\operatorname{\mathfrak{su}}}
\newcommand{\euu}{\operatorname{\mathfrak{u}}}

\newcommand{\w}{{\mathchoice{\,{\scriptstyle\wedge}\,}{{\scriptstyle\wedge}}
      {{\scriptscriptstyle\wedge}}{{\scriptscriptstyle\wedge}}}}
\newcommand{\lhk}{\mathbin{\hbox{\vrule height1.4pt width4pt depth-1pt 
             \vrule height4pt width0.4pt depth-1pt}}}

\newcommand{\be}{\begin{equation}}
\newcommand{\ee}{\end{equation}}
\newcommand{\bpm}{\begin{pmatrix}}
\newcommand{\epm}{\end{pmatrix}}

\numberwithin{equation}{section}

\newtheorem{theorem}{Theorem}

\newtheorem{proposition}{Proposition}
\newtheorem{corollary}{Corollary}

\theoremstyle{remark}
\newtheorem{definition}{Definition}

\newtheorem{remark}{Remark}

\begin{document}

\author[R. Bryant]{Robert L. Bryant}
\address{Duke University Mathematics Department\\
         P.O. Box 90320\\
         Durham, NC 27708-0320}
\email{\href{mailto:bryant@math.duke.edu}{bryant@math.duke.edu}}
\urladdr{\href{http://www.math.duke.edu/~bryant}%
         {http://www.math.duke.edu/\lower3pt\hbox{\symbol{'176}}bryant}}

\dedicatory{This article is dedicated to the
memory of Shiing-Shen Chern, whose beautiful works 
and gentle encouragement have had the most profound 
influence on my own research.}

\title[Almost complex 6-manifolds]
      {On the geometry\\
       of almost complex 6-manifolds}

\date{September 15, 2005}

\begin{abstract}

This article discusses some basic geometry 
of almost complex $6$-manifolds.

A $2$-parameter family 
of intrinsic first-order functionals 
on almost complex structures on $6$-manifolds
is introduced and their Euler-Lagrange equations 
are computed.

A natural generalization of holomorphic bundles 
over complex manifolds to the almost complex case
is introduced.  The general almost complex manifold
will not admit any nontrivial bundles of this type, but 
there is a class of nonintegrable almost complex
manifolds for which there are such nontrivial bundles.
For example, the $\Gtwo$-invariant almost complex
structure on the $6$-sphere admits such nontrivial 
bundles.  This class of almost complex manifolds 
in dimension~$6$ will be referred to as \emph{quasi-integrable}
and a corresponding condition for unitary structures
is considered.

Some of the properties of quasi-integrable structures
(both almost complex and unitary) are developed
and some examples are given.

However, it turns out that quasi-integrability 
is not an involutive condition, so the full 
generality of these structures in Cartan's sense
is not well-understood.  The failure of this
involutivity is discussed and some constructions
are made to show, at least partially, how general
these structures can be.
\end{abstract}



\thanks{
Thanks to Duke University for its support via
a research grant and to the NSF for its support
via grants DMS-9870164 and DMS-0103884.  
\hfill\break
\hspace*{\parindent} 
This is Version~$2.0$.  Version~$1$ was posted 
to the arXiv on August 23, 2005 and appears 
as {\tt math.DG/0508428}.
}

\maketitle

\setcounter{tocdepth}{2}
\tableofcontents

\section{Introduction}\label{sec: intro}
In the theory of almost complex manifolds,
dimension~$6$ is special for a number of reasons. 

First, there is a famous non-integrable almost complex structure 
on the $6$-sphere and the tantalizing problem of determining
whether or not there exists an integrable complex structure
on the $6$-sphere remains open as of this writing, 
despite much interest and effort.

Second, there exist nontrivial functorial constructions 
of top-degree differential forms associated to an almost 
complex structure~$J$ on a $6$-manifold that are first-order 
in the almost complex structure, something that does not happen
in dimension~$4$.  These are reviewed in~\S\ref{sec: ac6mfds}.

In particular, there is a $2$-parameter family of
first-order $6$-forms invariantly associated to an almost
complex structure in dimension~$6$ and hence these can be
regarded as defining a $2$-parameter family of natural 
Lagrangians for $6$-dimensional almost complex structures.  
The Euler-Lagrange equations of these functionals are computed 
and some of their properties and some examples are discussed.
Integrable almost complex structures are
critical points for these functionals, but it turns out
that there are many nonintegrable almost complex structures
that are also critical points for all of these functionals.

Third, there is a natural generalization of the notion 
of holomorphic bundles to almost complex manifolds (this
is true in any dimension), one that is based on replacing
the holomorphic transition function definition with 
the structure of a connection whose curvature is of type~$(1,1)$.
This opens the way to discussing a generalization of 
the notion of Hermitian-Yang-Mills connections to almost
complex manifolds.  I explain this in~\S\ref{sec: pdoholbundles}.

For the general nonintegrable almost complex manifold, 
this notion turns out to be too rigid to admit
any nontrivial examples, but, in dimension~$6$, there is
a special class of almost complex manifolds, which I have
called \emph{quasi-integrable}, for which the overdetermined
system that defines these generalized Hermitian-Yang-Mills
connections turns out to be just as well-behaved, locally,
as the familiar system on an integrable complex manifold.
In particular, important examples, such as the standard
nonintegrable structure on the $6$-sphere, turn out to be
quasi-integrable (more generally, nearly K\"ahler almost
complex manifolds turn out to be quasi-integrable).%
\footnote{It was Gang Tian, in 1998, who pointed 
out to me that the (pseudo-)Hermitian-Yang-Mills connections 
on complex vector bundles over~$S^6$ 
would be of interest in understanding the singular behavior 
of certain generalized Yang-Mills connections 
associated to calibrated geometries.  
I would like to thank him for his suggestion.}

My original goal was to regard quasi-integrable almost
complex structures as a generalization of integrable ones
to which the methods of Gromov~\cite{MR0864505} might be applicable
for establishing existence on compact $6$-manifolds. 
Then one might be able to use the category of pseudo-holomorphic 
bundles on such structures as a replacement for the
usual geometry of holomorphic mappings and submanifolds.
One might even be able to use the functionals described
above to look for particularly good quasi-integrable
structures on $6$-manifolds.

Such a program raises the natural question of how general 
or `flexible' the quasi-integrable almost complex 
structures are.  Unfortunately, quasi-integrability 
for almost complex structures turns out not to be a 
simple matter of imposing a smooth set of first-order 
equations on an almost complex structure.  Instead,
the space of $1$-jets of quasi-integrable almost complex
structures is a singular subset of codimension~$8$
in the space of all $1$-jets of all almost complex
structures, the singularity being the worst exactly
along the locus of $1$-jets of integrable almost
complex structures.  Moreover, unless the first Chern 
class~$c_1(TM,J)$ is $3$-torsion (a condition that depends 
only on the homotopy class of the almost complex 
structure~$J$), the $1$-jet of a quasi-integrable almost complex 
structure~$J$ on a compact manifold~$M^6$ cannot 
avoid this singular locus.  Such behavior does not 
bode well for application of Gromov's methods.

Instead it turns out to be better to work with a closely related notion, 
that of quasi-integrable $\U(3)$-structures on $6$-manifolds
(and their hyperbolic analogs, quasi-integrable $\U(1,2)$-structures), 
since this condition \emph{is} defined by a smooth set of~$16$ 
first-order equations on the $\U(3)$-structure 
(respectively, $\U(1,2)$-structure). 

In fact, what I have called 
the class of quasi-integrable $\U(3)$-structures 
was already considered by Hervella and Vidal~\cite{MR0431008,MR0577582} 
and by Gray and Hervella~\cite{MR0581924}, who developed
a nomenclature that would describe these~$\U(3)$-structures
as being of type~$\cW_1{\oplus}\cW_3{\oplus}\cW_4$. 

The relation between the two notions of quasi-integrability
is that, away from the singular locus, a quasi-integrable~$J$
has a canonical reduction to a quasi-integrable unitary structure~$(J,\eta)$
(where~$\eta$ is a nonsingular form of type~$(1,1)$ with respect to~$J$)
while, for any quasi-integrable unitary structure~$(J,\eta)$,
the underlying almost complex structure~$J$ is quasi-integrable.
(However, examples show that when~$J$ is quasi-integrable,
there may not be a nonsingular $\eta$ defined across the singular locus
such that~$(J,\eta)$ is quasi-integrable.)

Unfortunately, it turns out that although quasi-integrable
unitary structures are well-behaved up to second order, 
the \textsc{pde} system that defines them is \emph{not} involutive.
It is not even formally integrable, as obstructions appear
at the third order.  This was something that I did not know
(or even suspect) at the time of my lectures 
on this material in 1998 and 2000.
In fact, in my lecture in 2000, I erroneously
claimed that quasi-integrability for~$\U(3)$-structures was
an involutive condition.  It was only later, when I was 
writing up a careful proof up for publication, that I discovered
my mistake.  This is explained fully in~\S\ref{ssec: generality},
where a third-order obstruction to formal integrability is derived.

This third-order obstruction shows 
that the quasi-integrable~$\U(3)$-structures 
on a connected~$M^6$ split into two disjoint classes:  
Those for which the Nijenhuis tensor vanishes identically 
(and hence the underlying almost complex structure 
is integrable) and those for which the Nijenhuis tensor 
vanishes only along a proper subset of real codimension~$2$.
The former class is, of course, very well understood
and forms an involutive class of its own.  The latter class
is much less understood.  Indications are that the
latter class is, in fact, \emph{less} general than the
former class, as counter-intuitive as this may seem.%
\footnote{I realize that I am not being very precise
about the notion of `generality' being employed here.
Roughly speaking, what I mean is that, for $k\ge3$, 
the space of $k$-jets of structures in the former
class has larger dimension than the space of $k$-jets 
of structures in the latter class.}  For example, 
any smooth $1$-parameter deformation~$(J_t,\eta_t)$ 
through quasi-integrable~$\U(3)$-structures 
of the generic $\U(3)$-structure~$(J_0,\eta_0)$
with~$J_0$ integrable will necessarily have~$J_t$ be
integrable for $t$ sufficiently small.

I have not proved that the third-order obstruction that
I identify is the only third-order obstruction, much less
have I shown that, when this third-order obstruction is 
taken into account, there will be no higher order obstructions
to formal integrability.  This would require some rather
complex calculations that, so far, I have not been 
able to do.  This remains an interesting problem for the future.

\begin{remark}[Sources]
The results in this article were mostly reported on 
in two lectures~\cite{SemBesse1998,GrayConf2000} 
on aspects of the geometry of almost complex $6$-manifolds.
Of course, this written article contains many details 
and proofs that could not be given in those lectures for want of time.
Moreover, the results on noninvolutivity were not known
at the time of the lectures.

I discussed this work with Professor Chern 
from time to time and intended to write the results up for publication, 
but, after discovering the obstructions to formal integrabilty
discussed above, I intended to resolve that issue first.  
Unfortunately, I was not able to do so 
before moving on to other projects, 
so the write-up languished in my private files.  

Since I am no longer working in this area 
and interest continues to be expressed in my results, 
I have decided to release this article in its present form.

Those familiar with Professor Chern's unrivaled mastery 
of the method of equivalence and his use of differential forms
will see immediately what an enormous debt I owe to him
in this and many other works.  He will continue to be sorely
missed. 
\end{remark}

\begin{remark}[Verbitsky's eprint]
On 28 July 2005, Mikhail Verbitsky wrote to me and informed me of his
eprint~\cite{math.DG/0507179}, which contains several results that overlap
those of my own, though they were arrived at independently. 
The main overlaps appear to me to be these:  First, he 
rediscovered one of the two first-order functionals that 
I list in Theorem~\ref{thm: InvFuncts}.  Second, he
identified the class of almost complex structures that
I have called elliptically quasi-integrable as being of 
interest (though for different reasons than I have).  
Third, he showed that, among the strictly quasi-integrable 
almost complex structures, the only ones that are critical
for the functional that he defines are the ones that
underlie nearly K\"ahler structures.%
\footnote{His result is loosely related 
to my Proposition~\ref{prop: nrlypsdoKahler}, 
though our hypotheses and conclusions are quite different.}
\end{remark}

\section{Almost complex $6$-manifolds}
\label{sec: ac6mfds}
Let~$M^6$ be a $6$-manifold and let~$J:TM\to TM$ 
be an almost complex structure on~$M$.%
\footnote{Throughout this article, when~$M$ is a manifold, 
the space of complex-valued smooth alternating~$p$-forms 
on~$M$ will be denoted by~$\cA^p(M)$.} 
As usual, one has the splitting
\be
\cA^1(M) = \cA^{1,0}(M)\oplus\cA^{0,1}(M),
\ee
where~$\cA^{1,0}(M)$ consists of the complex-valued $1$-forms~$\alpha$
on~$M$ that satisfy~$\alpha(Jv) = \iC\,\alpha(v)$ for all~$v\in TM$
and~$\cA^{0,1}(M) = \overline{\cA^{1,0}(M)}$.  This induces
a type decomposition
\be
\cA^m(M) = \bigoplus_{p+q=m}\cA^{p,q}(M)
\ee
for~$m\ge0$ in the usual way, where
\be
\cA^{p,q}(M) = \Lambda^p\bigl(\cA^{1,0}(M)\bigr)
                \otimes \Lambda^q\bigl(\cA^{0,1}(M)\bigr).
\ee
(Strictly speaking, one should write~$\cA^{p,q}_J(M)$ for~$\cA^{p,q}(M)$, 
but this will be done only when the almost complex structure~$J$
is not clear from context.)

The exterior derivative~$\d:\cA^*(M)\to\cA^*(M)$
splits into a sum of terms~$\d^{r,s}:\cA^{p,q}\to\cA^{p+r,q+s}$
where, of course,~$r+s=1$ and, as is easy to see, $r,s\ge-1$.
Thus, one can write
\be
\d = \d^{2,-1} + \d^{1,0} + \d^{0,1} + \d^{-1,2}.
\ee
I will also use the common notation~$\d^{1,0}=\p$ and~$\d^{0,1}=\pbar$.

The Leibnitz rule implies that the 
operator~$\d^{-1,2}:\cA^{1,0}(M)\to\cA^{0,2}(M)$
is linear over the~$C^\infty$ functions on~$M$ and so 
represents a tensor.  This tensor is known as the
\emph{Nijenhuis} tensor and, by the famous
Newlander-Nirenberg Theorem~\cite{MR19577}, it vanishes identically if
and only if the almost complex structure is actually complex.

\begin{remark}[Existence]
In dimension~$6$, it is easy to determine when a given manifold~$M$
carries an almost complex structure.  The classical condition
in terms of the characteristic classes is that~$w_1(M) = \beta(w_2(M))=0$,
where~$w_1(M)$ and~$w_2(M)$ are the Stiefel-Whitney classes and
~$\beta:H^2(M,\bbZ_2)\to H^3(M,\bbZ)$ is the Bockstein homomorphism.
Equivalently, $M$ must be orientable (i.e., $w_1(M)=0$)
and~$w_2(M)$ must be of the form $w_2(M) \equiv c_1(K) \mod (2)$,
where~$K$ is some complex line bundle over~$M$, (i.e., $M$ must
possess a $\Spin^c(6)$-structure).

The set~$\cJ(M)$ of almost complex structures on~$M$ is the
space of sections of a bundle~$\Js(M)\to M$ whose fiber is modeled
on~$\GL(6,\bbR)/\GL(3,\bbC)$.  In fact, if~$\Fs(M)\to M$ is the 
right principal $\GL(6,\bbR)$-bundle of coframes~$u:T_xM\to\bbR^6$,
then~$\Js(M) = \Fs(M)/GL(3,\bbC)$.  Unless otherwise stated,
the topology on~$\cJ(M)$ will be the Whitney $C^\infty$ topology.

Each section~$J$ of~$\Js(M)$
determines an orientation of~$M$ and a line bundle~$K 
= \Lambda^3_\bbC(T^{1,0}_JM)^*$ such that~$w_2(M)\equiv c_1(K) \mod (2)$.
Moreover, when~$M$ is compact, the set~$\cJ(M,o,K)$ of sections of~$\Js(M)$ 
that determine a given orientation~$o$ and line bundle~$K$ 
is easily seen to be nonempty.
\end{remark}

\begin{remark}[First-order invariants]
\label{rem: 1jets}
The group~$\Diff(M)$ acts transitively on~$\Js(M)$ in the
obvious way, but it does not act transitively 
on~$J^1\bigl(\Js(M)\bigr)$.  

In fact, the $1$-jet of~$J$ at~$x\in M$ determines the 
complex linear mapping
\be
d^{-1,2}_x:(T^{1,0}_xM)^*\to \Lambda^{2,0}(T^{0,2}_xM)^*.
\ee
Conversely, given two almost complex structures~$I$ and~$J$
defined in a neighborhood of~$x$ 
that have the same $0$-jet at~$x$ (i.e., $I_x=J_x$) 
and that determine the same linear mapping~$d^{-1,2}_x$, 
there exists a diffeomorphism~$\phi:M\to M$ that is the identity at~$x$
to first-order and has the property that~$\phi^*(J)$ and~$I$ 
have the same $1$-jet at~$x$.

In this sense, the Nijenhuis tensor is a complete first-order
invariant of almost complex structures under the action of the 
diffeomorphism group.

It is, perhaps, also worth remarking that the map~$d^{-1,2}_x$
is essentially arbitrary.  (This is true in any even dimension, 
not just (real) dimension~$6$):  For any complex 
constants~$C^i_{\bj\bk}=-C^i_{\bk\bj}$, consider the complex
valued $1$-forms on~$\bbC^n$ defined by
\be
\alpha^i = \d z^i + {\ts\frac12}C^i_{\bj\bk}\,{\bar z}^j\,\d {\bar z}^k.
\ee
On an open neighborhood~$U\subset\bbC^n$ of~$z = 0$, 
these $n$ $1$-forms and their complex conjugates are linearly independent.  
As a consequence, there is a unique almost complex structure~$J_C$ on~$U$
for which the~$\alpha^i$ are a basis of the~$J_C$-linear $1$-forms on~$U$.
A straightforward calculation yields
\be
d^{-1,2}_0\bigl({\alpha^i}{\vphantom{p}\vrule}_0\bigr)
 = {\ts\frac12}C^i_{\bj\bk}\,\overline{\alpha^j}{\vphantom{p}\vrule}_0\w
          \overline{\alpha^k}{\vphantom{p}\vrule}_0\,,
\ee
thus verifying that~$d^{-1,2}_0$ can be arbitrarily prescribed.
\end{remark}

\subsection{Invariant forms constructed from the Nijenhuis tensor}
\label{ssec: Nijenhuisinvariants}
Let~$U\subset M$ be an open set on which there exist
linearly independent $1$-forms~$\alpha^1,\alpha^2,\alpha^3$
in~$\cA^{1,0}(U)$.  Write~$\alpha = (\alpha^i):TU\to \bbC^3$
and note that~$\alpha$ maps each fiber~$T_pU$ isomorphically
(and complex linearly) onto~$\bbC^3$.  I will say that~$\alpha$
is a $J$-complex local coframing on~$M$ with domain~$U$.

There exists a unique smooth mapping~$N(\alpha):U\to M_{3\times3}(\bbC)$
such that
\be
\label{eq: Nalphadefined}
\d^{-1,2}\alpha=\d^{-1,2}\bpm\alpha^1\\ \alpha^2\\ \alpha^3\epm
       = N(\alpha)\, 
    \bpm\overline{\alpha^2\w\alpha^3}\\ 
        \overline{\alpha^3\w\alpha^1}\\ 
        \overline{\alpha^1\w\alpha^2}\epm.
\ee
If~$\beta = (\beta^i):TU\to\bbC^3$ is any other $J$-linear
coframing on~$U$, then there exists a unique function~$g:U\to\GL(3,\bbC)$
such that~$\beta = g\alpha$ and computation yields
\be
\label{eq: Ntransition}
N(\beta) = N(g\alpha) = \det(\bar g)^{-1}\,g\,N(\alpha)\,{}^t\bar g\,.
\ee

This motivates investigating the representation~$\rho:\GL(3,\bbC)
\to \End\bigl(M_{3\times3}(\bbC)\bigr)$ defined by
\be
\rho(g)(a) = \det(\bar g)^{-1}\,g\,a\,{}^t\bar g.
\ee
In particular, it will be important to understand its orbit structure.%

\subsubsection{The form~$\Phi$}
Clearly, one has $\bigl|\det(\rho(g)a)\bigr|^2=|\det(g)|^{-2}|\det(a)|^2$,
so that
\be
\bigl|\det\bigl(N(g\alpha)\bigr)\bigr|^2
= |\det(g)|^{-2}\,\bigl|\det\bigl(N(\alpha)\bigr)\bigr|^2,
\ee
implying that the nonnegative $6$-form
\be
\Phi = {\ts\frac{\iC}{8}}\,\bigl|\det\bigl(N(\alpha)\bigr)\bigr|^2
\,\alpha^1\w\alpha^2\w\alpha^3\w\overline{\alpha^1}\w\overline{\alpha^2}
\w\overline{\alpha^3}
\ee
is well-defined on~$U$ independent of the choice of~$\alpha$.
Consequently,~$\Phi$ is the restriction to~$U$ of
a nonnegative $6$-form, also denoted~$\Phi$, that
is well-defined globally on~$M$.  
When it is necessary to make the dependence on the almost complex 
structure~$J$ explicit, I will denote this $6$-form by~$\Phi(J)$.

\subsubsection{The form~$\omega$}
Another invariant is the canonical $(1,1)$-form~$\omega$ 
that is defined as follows:  Let
\be
Q:M_{3\times3}(\bbC)\times M_{3\times3}(\bbC)
\to M_{3\times3}(\bbC)
\ee
be the unique bilinear map satisfying
$Q(a,b) = Q(b,a)$ and~$Q(a,a) = \det(a) a^{-1}$.  
This~$Q$ satisfies~$Q({}^ta,{}^tb) = {}^tQ(a,b)$ and
\be
Q\bigl(gah,gbh\bigr)
= \det(g)\det(h)\,h^{-1}\,Q(a,b)\,g^{-1}
\ee
for all~$g,h\in\GL(3,\bbC)$ and~$a,b\in M_{3\times3}(\bbC)$.
As a result, the $(1,1)$-form~$\omega$ defined on~$U$ by
\be
\label{eq: omegadefined}
\omega = {\ts\frac{\iC}2}\,\,{}^t\alpha\w
        Q\bigl({}^t\!N(\alpha),\overline{N(\alpha)}\,\bigr)\w\overline{\alpha},
\ee
does not depend on the choice of the coframing~$\alpha$ and,
since~${}^tQ({}^ta,\overline a) = \overline{Q({}^ta,\overline a)}$,
this form is real-valued.  
Consequently, it is the restriction to~$U$
of a real-valued $(1,1)$-form that is well-defined globally on~$M$.
When it is necessary to make the dependence on~$J$ explicit,
this $2$-form will be denoted by~$\omega(J)$.

\begin{remark}[Other dimensions]
\label{rem: otherdims}
Although this will have no bearing on the rest of this
article, it should be remarked that a corresponding construction
of a canonical $(1,1)$-form exists in other dimensions as well:
If~$J$ is an almost complex structure on~$M^{2n}$, then,
for any local $J$-linear coframing~$\alpha = (\alpha^i):TU\to\bbC^n$,
one will have
\be
\d^{-1,2}\alpha^i = {\ts\frac12} N^i_{\bj\bk}\,
                      \overline{\alpha^j}\w\overline{\alpha^k}
\ee
where~$N^i_{\bj\bk}=-N^i_{\bk\bj}$.  Then the $(1,1)$-form
\be
\omega = -{\ts\frac\iC4} N^i_{\bar\ell\bk}\overline{N^\ell_{\bi\bj}}\,\,
               \alpha^j\w\overline{\alpha^k}
\ee
is easily seen to be real-valued and well-defined, independent of
the choice of $\alpha$.  When~$n=3$, this reduces to the
definition of~$\omega$ already given above.
\end{remark}

\subsubsection{The form~$\psi$}
A third invariant can be defined as follows:  Let
\be
P:M_{3\times3}(\bbC)\times M_{3\times3}(\bbC)\times M_{3\times3}(\bbC)\to\bbC
\ee
be the symmetric trilinear map satisfying~$P(a,a,a)=\det(a)$
for~$a\in M_{3\times3}(\bbC)$.  
Note the identities~$P({}^ta,{}^tb,{}^tc)=P(a,b,c)$ and
\be
P(a,b,c) = {\ts\frac13}\tr\bigl(Q(a,b)c\bigr) 
         = {\ts\frac13}\tr\bigl(Q(a,c)b\bigr),
\ee
which will be useful below. The map~$P$ has the equivariance
\be
P(gah,gbh,gch) = \det(g)\det(h)\,P(a,b,c)
\ee
for all~$g,h\in\GL(3,\bbC)$.
Thus, the $(3,0)$-form on~$U$ defined by
\be
\psi = P\bigl(\,\overline{N(\alpha)},\overline{N(\alpha)},
             {}^t\!N(\alpha)\,\bigr)
        \,\alpha^1\w\alpha^2\w\alpha^3,
\ee
is independent of the choice of coframing~$\alpha$ on~$U$.  Consequently,
it is the restriction to~$U$ of a well-defined $(3,0)$-form on~$M$.  
When it is necessary to make the dependence on~$J$ explicit, 
this $(3,0)$-form will be denoted by~$\psi(J)$.  

\subsubsection{Tautological forms}
One way to interpret the above constructions 
is that there exist, on~$J^1\bigl(\Js(M)\bigr)$, 
unique, $M$-semi-basic, smooth differential forms~$\Phib$ 
(of degree~$6$), $\psib$ (of degree~$3$), and~$\omegab$ (of degree~$2$)  
such that, for any almost complex structure~$J:M\to\Js(M)$, 
one has
\be
\label{eq: tautonJ1J}
j^1(J)^*(\Phib)=\Phi(J),\qquad
j^1(J)^*(\psib)=\psi(J),\qquad
j^1(J)^*(\omegab)=\omega(J),
\ee
where $j^1(J):M\to J^1\bigl(\Js(M)\bigr)$ 
is the $1$-jet lifting of~$J:M\to\Js(M)$.

\subsubsection{Bi-forms}
For later use, I will point out the existence
of two further first-order invariant `bi-form' tensors that will 
turn out to be important. These two forms~$E_1(J)$ and~$E_2(J)$ 
are sections of the bundle~$A^{1,1}(M)\otimes A^{3,0}(M)$ and
are given in terms of a local complex coframing~$\alpha$ by the
formulae
\be
\label{eq: E1J}
E_1(J) 
= {\ts\frac\iC2}
  \,{}^t\alpha\w Q\bigl({}^t\!N(\alpha),{}^t\!N(\alpha)\,\bigr)\w\overline{\alpha}
        \otimes P\bigl(\,\overline{N(\alpha)},\overline{N(\alpha)},
          \overline{N(\alpha)}\,\bigr)\,\alpha^1\w\alpha^2\w\alpha^3
\ee
and
\be
\label{eq: E2J}
E_2(J) 
= {\ts\frac\iC2}
 \,{}^t\alpha\w Q\bigl(\overline{N(\alpha)},\overline{N(\alpha)}\,\bigr)
             \w\overline{\alpha}
        \otimes P\bigl(\,{}^t\!N(\alpha),{}^t\!N(\alpha),
          \overline{N(\alpha)}\,\bigr)\,\alpha^1\w\alpha^2\w\alpha^3.
\ee
The change-of-coframing formulae already discussed show that these
bi-forms are globally defined on~$M$.  
The reader will note the similarity of these definitions to
\be
\label{eq: E0J}
\begin{aligned}
E_0(J) &= \omega(J)\otimes\psi(J),\\
&= {\ts\frac\iC2}
 \,{}^t\alpha\w Q\bigl(\,{}^t\!N(\alpha),\overline{N(\alpha)}\,\bigr)
             \w\overline{\alpha}
        \otimes P\bigl(\overline{N(\alpha)},\overline{N(\alpha)},
          {}^t\!N(\alpha)\bigr)\,\alpha^1\w\alpha^2\w\alpha^3,
\end{aligned}
\ee
which is also a section of~$A^{1,1}(M)\otimes A^{3,0}(M)$.

\subsection{Relations and Inequalities}
\label{ssec: ineqs}
There are some relations among these invariant forms:

\begin{proposition}[Relations]
\label{prop: invineq}
The invariants~$\omega$, $\psi$, and~$\Phi$ satisfy 
\be
\label{eq: Psiformula}
{\ts\frac98}\iC\,\psi\w\overline{\psi} = \Phi + {\ts\frac43}\omega^3.
\ee
In particular,
\be
\label{eq: omega3Phiineq}
-{\ts\frac18}\,\Phi\le {\ts\frac16}\,\omega^3.
\ee
\end{proposition}

\begin{proof}
The equation \eqref{eq: Psiformula} is equivalent to 
\be
\label{eq: PdetQreln}
9\bigl|P(a,a,{}^t\bar a)\bigr|^2 
= \bigl|\det(a)\bigr|^2 + 8\,\det\bigl(Q(a,{}^t{\bar a})\bigr)
\ee
for~$a\in M_{3\times3}(\bbC)$
while the inequality~\eqref{eq: omega3Phiineq} 
is equivalent to the homogeneous polynomial inequality
\be
\label{eq: detQineq}
-{\ts\frac18}\,\bigl|\det(a)\bigr|^2 \le \det\bigl(Q(a,{}^t{\bar a})\bigr),
\ee
Obviously~\eqref{eq: PdetQreln}
implies~\eqref{eq: detQineq}, 
so it suffices to establish~\eqref{eq: PdetQreln}.
This latter relation can be proved directly simply by expanding
out both sides, but this is something of a mess and is somewhat
unconvincing.  

To establish~\eqref{eq: PdetQreln}, it suffices to prove it 
on an open subset of~$M_{3\times3}(\bbC)\simeq\bbC^9$.  
Moreover, because of the 
equivariance of the mappings~$P$ and~$Q$, if this relation
holds for a given~$a\in M_{3\times3}(\bbC)$, then it holds
on the entire~$\rho$-orbit of~$a$.  Thus, I will describe a 
`normal form' for the $\rho$-action on an open set 
in~$M_{3\times3}(\bbC)$.  (This normal form will also come in 
handy for other reasons.)

First, consider the restriction of the~$\rho$-action 
to~$\SL(3,\bbC)\subset\GL(3,\bbC)$. This restricted action 
is just~$\rho(g)(a) = g a {}^t\bar g$ and hence~$M_{3\times3}(\bbC)$ 
is reducible as a real $\SL(3,\bbC)$-representation.  
Namely, write~$a = h + \iC k$ where~$h$ and~$k$ 
are Hermitian symmetric, $h = {}^t\bar h$ and $k = {}^t\bar k$, 
and note that this splitting is preserved by
the $\rho$-action of~$\SL(3,\bbC)$.

As is well-known in linear algebra, there is an open set
in the space of pairs of $3$-by-$3$ Hermitian symmetric matrices 
that can be simultaneously diagonalized by an element of~$\SL(3,\bbC)$.%
\footnote{For example, any pair~$(h,k)$ for which the cubic polynomial
$p(t) = \det(h + t k)$ has three real distinct roots can be simultaneously
diagonalized by this action and the set of such pairs is open.}  
In particular, a nonempty open set of elements of~$M_{3\times3}(\bbC)$ 
are $\rho$-equivalent to diagonal elements of the form
\be
a = \bpm a_1&0&0\\
          0&a_2&0\\ 
          0&0&a_3\\  \epm
\ee
where the~$a_i$ are nonzero.  Now, after $\rho$-acting by
a diagonal element of~$\GL(3,\bbC)$, one can clearly reduce such
an element to the form
\be
\label{eq: anormalized}
a = \bpm \E^{\iC\lambda_1}&0&0\\
          0&\E^{\iC\lambda_2}&0\\ 
          0&0&\E^{\iC\lambda_3}\\  \epm
\ee
where~$\lambda_1+\lambda_2+\lambda_3 = 0$.  Thus, the union of the
$\rho$-orbits of such elements contains an open set in~$M_{3\times3}(\bbC)$, 
so it suffices to prove~\eqref{eq: PdetQreln} for~$a$ of this form.   

However, for~$a$ in the form~\eqref{eq: anormalized},
one has~$\det(a)=1$ and computation yields
\be
P(a,a,{}^t\bar a) = {\ts\frac13}\bigl(
\E^{-2\iC\lambda_1}+\E^{-2\iC\lambda_2}+\E^{-2\iC\lambda_3}\bigr)
\ee
and
\be
Q(a,{}^t\bar a) = \bpm \cos(\lambda_2{-}\lambda_3)&0&0\\
          0&\cos(\lambda_3{-}\lambda_1)&0\\ 
          0&0&\cos(\lambda_1{-}\lambda_2)\\  \epm.
\ee
The relation~\eqref{eq: PdetQreln} for~$a$ 
of the form~\eqref{eq: anormalized} now follows trivially.
\end{proof}

\begin{remark}[More identities]
Although the structure equations will be explored more thoroughly
in a later section, the reader may be interested to note the
identity
\be
\d^{-1,2}\omega = \frac{\iC}2\,\tr\bigl({}^tN(\alpha)Q({}^tN(\alpha),
                    \overline{N(\alpha)}\,\bigr)\,
                \overline{\alpha^1\w\alpha^2\w\alpha^3},
\ee
which follows immediately from~\eqref{eq: Nalphadefined} 
and~\eqref{eq: omegadefined}.

Since the identity~$\tr\bigl({}^ta\,Q({}^ta,\bar a)\bigr) 
= 3 P(a,a,{}^t\bar a)$ holds for all~$a\in M_{3\times 3}(\bbC)$,
one has
\be
\label{eq: dm12omega}
\d^{-1,2}\omega = {\ts\frac32}\iC\,\overline{\psi},
\ee
so that one has the useful identity
\be
\label{eq: domega}
\d\omega = 3\Im(\psi) + \p\omega + \pbar\omega 
         = 3\Im(\psi) + 2\Re(\pbar\omega).
\ee
Of course,~$\pbar\omega=\d^{0,1}\omega =\overline{\p\omega}$ 
is a second-order invariant of~$J$.

Similarly, the $(2,2)$-form~$\d^{-1,2}\psi$ is first-order in~$J$ 
while the expression~$\pbar\psi$ is second-order.  
In fact, with respect to any complex coframing~$\alpha$, one has
\be
\label{eq: dm12psi}
\d^{-1,2}\psi 
= \bpm\alpha^2\w\alpha^3&\alpha^3\w\alpha^1&\alpha^1\w\alpha^2\epm
R\bigl(N(\alpha)\bigr)
\bpm\overline{\alpha^2\w\alpha^3}\strut\\
    \overline{\alpha^3\w\alpha^1}\strut\\
    \overline{\alpha^1\w\alpha^2}\strut
\epm
\ee
where~$R:M_{3\times3}(\bbC)\to M_{3\times3}(\bbC)$ is the function
\be
\label{eq: Rdefined}
R(a) = P(\bar a,\bar a,{}^ta)\,a.
\ee
Thus, the identity~$\tr\bigl({}^ta\,Q({}^ta,\bar a)\bigr) 
= 3P(a,a,{}^t\bar a)$ implies 
\be
\omega \w \d^{-1,2}\psi = {\ts\frac{3}2}\iC\,\psi\w\overline\psi,
\ee
which also follows from~\eqref{eq: dm12omega}
by differentiating the identity~$\omega\w\psi = 0$.

It seems likely (though I have not written out a proof) that the
ring of smooth first-order invariant forms definable for an almost complex
structure~$J$ on a $6$-manifold is generated by~$\omega$, 
$\psi$,~$d^{-1,2}\psi$, $\Phi$, and their conjugates.  

Note, by the way, that~\eqref{eq: Psiformula} 
implies that all of these first-order forms can be generated from~$\omega$
by repeatedly applying the operator~$d^{-1,2}$, conjugation, 
and exterior multiplication.  Thus, in some sense, 
$\omega$ is the fundamental object.%
\footnote{Of course, in other dimensions as well, the $(1,1)$-form~$\omega$ 
can be used to generate a collection of invariant first-order
forms.  See Remark~\ref{rem: otherdims}.}
\end{remark}

\begin{remark}[Nonvanishing invariants]
If $J$ is an integrable almost complex structure, then~$d^{-1,2}=0$,
so, of course, all of the invariant forms constructed above from the 
Nijenhuis tensor vanish.  In the other direction, one might be interested 
to know when it is possible to have one of these invariant 
forms be nowhere vanishing.  

Now, $\Phi(J)$ is nowhere vanishing if and only if the bundle map
\be
\d^{-1,2}:(T^{1,0}M)^*\to \Lambda^2\bigl(T^{0,1}M\bigr)^*
\ee
is an isomorphism.  This would imply that the bundles
$(T^{1,0}M)^*$ and~$\Lambda^2\bigl(T^{0,1}M\bigr)^*$ have
the same Chern classes, i.e., that $3c_1(J)=0$
and~$c_1(J)^2 = c_1(J)c_2(J)=0$ (where, of course, the Chern
classes are computed for the bundle~$TM$ with the complex structure~$J$).

Thus, for example, any almost complex structure on~$\bbC\bbP^3$
that is homotopic to the standard complex structure cannot have
a nonvanishing~$\Phi$, since the first Chern class
of the standard almost complex structure is nontrivial.

If~$\psi(J)$ were to be nonvanishing, 
then one clearly would have~$c_1(J)=0$, 
for, in this case, the `canonical bundle'~$K = \Lambda^3(T^{1,0}M)^*$
would be trivial (and not just $3$-torsion) because~$\psi(J)$
would be a nonvanishing section of it.

Finally, note that, if~$\omega(J)$ were everywhere nondegenerate,
then either~$\omega(J)^3>0$, in which case, by~\eqref{eq: Psiformula}
one sees that~$\psi(J)$ is nowhere vanishing,
or else~$\omega(J)^3<0$, in which case, again by~\eqref{eq: Psiformula},
one sees that~$\Phi(J)$ is nowhere vanishing.

One has a stronger result in the case that~$\omega(J)$
is a positive $(1,1)$-form:  It is easy to prove that, if~$Q({}^t\!a,\bar a)$
is a positive definite Hermitian symmetric matrix, then~$a$
is $\rho$-equivalent to a matrix of the form~\eqref{eq: anormalized}
in which all of the~$\lambda_i$ satisfy~$|\lambda_i|<\frac12\pi$.  
Consequently, if~$\omega(J)$ is positive, then both~$\psi(J)$
and~$\Phi(J)$ are nonvanishing and one has the inequality%
\footnote{N.B.:  The assumption that $\omega(J)$ be a positive $(1,1)$-form
is essential for the inequality~\eqref{eq: omposupperbound},
since it does not hold in general.}
\be
\label{eq: omposupperbound}
{\ts\frac16}\omega(J)^3 \le \Phi(J).
\ee
The case of equality in this last inequality will turn out to be
important below, where it will be seen to be equivalent to 
\emph{strict quasi-integrability} (when~$\omega$ is positive).
\end{remark}

\subsection{First-order functionals}
\label{ssec: firstorderfuncts}
The above constructions yield the following result:

\begin{theorem}[Invariant functionals]
\label{thm: InvFuncts}
If~$M^6$ is compact and oriented, then for any 
constants~$c_1$ and $c_2$, the functional
\be
\label{eq: Fcdefined}
F_c(J) = \int_M c_1\,\Phi(J) 
         + c_2\,{\ts\frac\iC8}\,\psi(J)\w\overline{\psi(J)}
\ee
is a diffeomorphism-invariant, first-order functional 
on the space~$\cJ_+(M)$ of almost complex structures on~$M$ 
that induce the given orientation.  
The functional~$F_c$ belongs to the class~$L^6_1$. \qed
\end{theorem}

\begin{remark}[A classification]
It is not difficult to show that any smooth function~$L:M_{3\times3}(\bbC)
\to\bbR$ that satisfies
\be
L\bigl(\rho(g)a\bigr) = |\det(g)|^{-2} L(a)
\ee
for all~$g\in \GL(3,\bbC)$ is of the form
\be
L(a) = c_1\,\bigl|\det(a)\bigr|^2 + c_2\,\bigl|P(\bar a,\bar a,{}^ta)\bigr|^2
\ee
for some real constants~$c_1$ and $c_2$.  From this and the
discussion in~Remark~\ref{rem: 1jets},
it follows without difficulty that the functionals listed 
in~Theorem~\ref{thm: InvFuncts} are the only smooth, 
diffeomorphism-invariant, first-order functionals on~$\cJ_+(M)$.

Note that the map~$(x,y)= \bigl(|\det(a)|^2,
|P(\bar a,\bar a,{}^ta)|^2\bigr):M_{3\times3}(\bbC)\to\bbR^2$
has, as its image, the entire closed quadrant~$x\ge0,\,y\ge0$, 
so that neither of the terms in~$L$ dominates the other in general.

The fact that these first-order functionals are in~$L^6_1$ 
should be no surprise;  as is usual with geometrically interesting 
functionals, they are on the boundary of Sobolev embedding.  
\end{remark}

\subsection{The Euler-Lagrange system}
It is natural to ask about the critical points of the functionals~$F_c$
introduced in Theorem~\ref{thm: InvFuncts}.  Because these functionals
are first-order, their Euler-Lagrange equations can be computed
using the method of Poincar\'e-Cartan forms~\cite{MR1985469}.

	To implement this method, I will first exhibit
the $1$-jet bundle~$J^1\bigl(\Js(M)\bigr)\to M$ 
as a quotient of a prolonged coframe bundle~$\Fs^{(1)}(M)\to M$ 
and then use the canonical forms on~$\Fs^{(1)}(M)$ to do the 
necessary computations.

\subsubsection{$(0,2)$-connections}
To motivate this, I will begin by reminding the reader about 
$(0,2)$-connections associated to almost complex structures.
If~$J$ is an almost complex structure on~$M^6$, then one has
the $\GL(3,\bbC)$-bundle~$\Fs(J)$ of~$J$-linear coframes~$u:T_xM\to\bbC^3$
and the canonical (i.e., tautological) $\bbC^3$-valued~$1$-form~$\eta$
on~$\Fs(J)$.  A connection on~$\Fs(J)$ is a $\GL(3,\bbC)$-equivariant
$1$-form~$\kappa$ on~$\Fs(J)$ with values in~$\eugl(3,\bbC)$ such
that the $2$-form~$\d\eta + \kappa\w\eta$ is semi-basic for the 
projection~$\Fs(J)\to M$.  One says that such a connection is 
a~\emph{$(0,2)$-connection} if it satisfies the stronger condition
that
\be
\d\eta^i = -\kappa^i_j\w\eta^j 
             + {\ts\frac12}C^i_{\bj\bk}\,\overline{\eta^j\w\eta^k}.
\ee
for some (necessarily unique) functions~$C^i_{\bj\bk} = - C^i_{\bk\bj}$
on~$\Fs(J)$. Such connections always exist and, in fact, the set of such 
connections is the space of sections 
of an affine bundle~$\Fs^{(1)}(J)\to\Fs(J)$ of real rank~$36$.

\subsubsection{The prolonged coframe bundle}
For the next step, recall that~$\GL(6,\bbR)$
can be embedded into~$\GL(6,\bbC)$ as the subgroup of invertible matrices
of the form
\be
\label{eq: GL6RinGL6C}
\bpm A & B \\ \bar B & \bar A\epm
\ee
where~$A$ and~$B$ lie in~$M_{3\times 3}(\bbC)$.  

Now, let~$\pi:\Fs(M)\to M$ 
be the bundle of~$\bbC^3$-valued coframes $u:T_xM\to\bbC^3$.
When there is no danger of confusion, I will simply write~$\Fs$ for~$\Fs(M)$.
Let~$\omega:T\Fs\to\bbC^3$ be the canonical $\bbC^3$-valued $1$-form
on~$\Fs$.  Sometimes, $\omega=(\omega^i)$ is called the `soldering form'.

Using the identification of~$\eugl(6,\bbR)$ with pairs of matrices
in~$\eugl(3,\bbC)$ implicit in~\eqref{eq: GL6RinGL6C}, one can see
that specifying a pseudo-connection on~$M$ is equivalent to specifying
a pair of $1$-forms~$\alpha$ and~$\beta$ 
with values in~$M_{3\times3}(\bbC)$ such that the $2$-form
\be
\d\omega + \alpha\w\omega + \beta\w\bar\omega
\ee
is semi-basic.  I will say that such a pseudo-connection is 
\emph{$\bbC$-compatible} if it satisfies the stronger condition
that
\be
\label{eq: Ccompdefined}
\d\omega^i + \alpha^i_j\w\omega^j + \beta^i_{\bj}\w\overline{\omega^j}
= {\ts\frac12} T^i_{\bj\bk}\,\overline{\omega^j\w\omega^k}
\ee
for some functions~$T^i_{\bj\bk}=-T^i_{\bk\bj}$ on~$\Fs$.  (Note
that it is necessary to consider pseudo-connections here because
the only actual connections satisfying such a condition are the
ones that satisfy~$T^i_{\bj\bk}=0$.)

As is easy to see, the pseudo-connections on~$\Fs$ 
that are $\bbC$-compatible are the sections of an affine bundle
~$\Fs^{(1)}\to\Fs$ of real rank~$144$.   On~$\Fs^{(1)}$, there
exist canonical complex-valued~$1$-forms~$\phi^i_j$ and~$\theta^i_{\bj}$
and complex-valued functions~$N^i_{\bj\bk}=-N^i_{\bk\bj}$ such
that the following structure equations hold
\be
\label{eq: 1ststreqonF1}
\d\omega^i = -\phi^i_j\w\omega^j - \theta^i_{\bj}\w\overline{\omega^j}
              + {\ts\frac12} N^i_{\bj\bk}\,\overline{\omega^j\w\omega^k}.
\ee
These canonical forms and functions are defined by the property that, 
when a $\bbC$-compatible pseudo-connection defined by an~$\alpha$ 
and~$\beta$ as in~\eqref{eq: Ccompdefined} 
is regarded as a section~$\nabla:\Fs\to\Fs^{(1)}$, then~$\nabla$ pulls
back~$\phi$ to be~$\alpha$, $\theta$ to be~$\beta$, and~$N$ to be~$T$.

It is convenient to think of~$N$ as taking values in $M_{3\times 3}(\bbC)$
as earlier, i.e.,~$N:\Fs^{(1)}\to M_{3\times3}(\bbC)$, since~$\GL(3,\bbC)$
clearly acts on~$\Fs^{(1)}$ in such a manner that this mapping is 
$\rho$-equivariant.

Given a choice of $\bbC$-compatible pseudo-connection on~$\Fs$ defined
by $1$-forms~$\alpha$ and~$\beta$ as above (and, therefore, defining
the function~$T$ on~$\Fs$), one can consider the plane field~$\beta^\perp
\subset T\Fs$ of corank~$18$ on~$\Fs$ defined by the 
equations~$\beta=(\beta^i_{\bj})=0$.  At~$u\in\Fs$, the plane~$\beta^\perp_u
\subset T_u\Fs$ is tangent to the $\GL(3,\bbC)$-structure associated
to an almost complex structure.  In fact, $\beta^\perp_u$ determines
a $1$-jet of an almost complex structure at~$u$.  Thus, there 
is a natural submersion~$\xi:\Fs^{(1)}\to J^1\bigl(\Js(M)\bigr)$.
The map~$\xi$ has been defined in such a way that it pulls back the
canonical contact system on~$J^1\bigl(\Js(M)\bigr)$ to be the ideal~$\cI$
generated by the $1$-forms~$\theta^i_{\bj}$.  In fact, the maximal integral 
manifolds of~$\cI$ on which
\be
\Omega = {\ts\frac\iC8}\,\omega^1\w\omega^2\w\omega^3\w
           \overline{\omega^1\w\omega^2\w\omega^3}\not=0,
\ee
are precisely the natural embeddings of the bundles~$\Fs^{(1)}(J)$
into~$\Fs^{(1)}$ for almost complex structures~$J$.

Now computations that could be carried out abstractly 
on~$J^1\bigl(\Js(M)\bigr)$ can be carried out explicitly on~$\Fs^{(1)}$
using the canonical forms and functions.  It is this technique that
will now be employed.

Recall that there are canonical forms~$\Phib$ and~$\psib$ 
on~$J^1\bigl(\Js(M)\bigr)$ that satisfy~\eqref{eq: tautonJ1J}.  
These forms pull back to~$\Fs^{(1)}$ via~$\xi$ to satisfy
\be
\xi^*\Phib = |\det(N)|^2\,\Omega
\ee
and
\be
\xi^*\psib = 
        P(\overline{N},\overline{N},{}^tN)\,\omega^1\w\omega^2\w\omega^3.
\ee
In particular, for any constants~$c_1$ and~$c_2$, one has
\be
\xi^*(c_1\Phib+c_2{\ts\frac\iC8}\,\psib\w\overline\psib) 
=\bigl(c_1|\det(N)|^2+c_2|P(\overline{N},\overline{N},{}^tN)|^2\bigr)\,\Omega.
\ee

Set~$L = L(N) = c_1 |\det(N)|^2+c_2|P(\overline{N},\overline{N},{}^tN)|^2$
and note that~$L$ has the equivariance
\be
L\bigl(\rho(g)a\bigr) = |\det(g)|^{-2}L(a)
\ee
for all~$a\in M_{3\times3}(\bbC)$.  I am now going to compute the
Euler-Lagrange system for the functional
\be
\Lambda = L(N)\,\Omega.
\ee

To do this effectively, I will need the second structure equations 
on~$\Fs^{(1)}$.  To compute these, write
\be
\d\omega^i
 = -\phi^i_j\w\omega^j 
      - \bigl(\theta^i_{\bj}+{\ts\frac12}
            N^i_{\bj\bk}\,\overline{\omega^k}\bigr)\w\overline{\omega^j}.
\ee
and, temporarily, set~$\psi^i_{\bj} 
= \theta^i_{\bj}+{\ts\frac12}N^i_{\bj\bk}\,\overline{\omega^k}$, so 
that the above equation takes the form
\be
\d\omega^i = -\phi^i_j\w\omega^j - \psi^i_{\bj}\w\overline{\omega^j}.
\ee
Taking the exterior derivative of this expression gives
\be
0 = -(\d\phi^i_j+\phi^i_k\w\phi^k_j+\psi^i_\bk\w\overline{\psi^k_\bj})
            \w\omega^j
    -(\d\psi^i_\bj + \phi^i_k\w\psi^k_\bj + \psi^i_{\bk}\w\overline{\phi^k_j})
           \w\overline{\omega^j},
\ee
so that, by Cartan's Lemma, there must exist $1$-forms~$\pi^i_{jk}=\pi^i_{kj}$,
$\pi^i_{j\bk}$, and~$\pi^i_{\bj\bk}=\pi^i_{\bk\bj}$ so that
\be
\begin{aligned}
\d\phi^i_j+\phi^i_k\w\phi^k_j+\psi^i_\bk\w\overline{\psi^k_\bj}
 & = -\pi^i_{jk}\w\omega^k - \pi^i_{j\bk}\w\overline{\omega^k},\\
\d\psi^i_\bj + \phi^i_k\w\psi^k_\bj + \psi^i_{\bk}\w\overline{\phi^k_j}
 & = - \pi^i_{k\bj}\w\omega^k -\pi^i_{\bj\bk}\w\overline{\omega^k}.
\end{aligned}
\ee
(The forms~$\pi$ are \emph{not} canonically defined on~$\Fs^{(1)}$, 
but this ambiguity will not cause a problem.)  Setting
\be
\nu^i_{\bj\bk} = \d N^i_{\bj\bk} + N^\ell_{\bj\bk}\,\phi^i_\ell
                 - N^i_{\bar\ell\bk}\,\overline{\phi^\ell_j}
                 - N^i_{\bj\bar\ell}\,\overline{\phi^\ell_k}
 = - \nu^i_{\bk\bj}\,,
\ee
one finds that the above equations imply
\be
\d\theta^i_\bj
\equiv 
-\phi^i_k\w\theta^k_\bj - \theta^i_\bk\w\overline{\phi^k_j}    
-\bigl(\pi^i_{\bj\bk}+{\ts\frac12}\nu^i_{\bj\bk}\bigr)\w\overline{\omega^k}
\mod \omega^1,\omega^2,\omega^3.
\ee

\subsubsection{The Euler-Lagrange ideal}
Now, given~$L= L(N) = \bar L$ as above, 
let~$L^{\bj\bk}_i = - L^{\bk\bj}_i$ be the functions such that
\be
\label{eq: Ljkidefined}
\d L = \Re\bigl(L^{\bj\bk}_i\,\d N^i_{\bj\bk}\bigr)
\ee
(Since~$L$ is a polynomial of degree~$6$ on~$M_{3\times3}(\bbC)$,
the functions~$L^{\bj\bk}_i$ are polynomials of degree~$5$.)

The equivariance of~$L$ and~$\Omega$ and the definitions
made so far imply that
\be
\d\bigl(L\,\Omega\bigr)
= \Re\bigl(L^{\bj\bk}_i\,\nu^i_{\bj\bk}\bigr)\w\Omega.
\ee

Let~$\omega_{[\bi]}$ be the $5$-form that satisfies
\be
\overline{\omega^j}\w\omega_{[\bi]} = \delta^j_i\,\Omega
\ee
and~$\omega^j\w\omega_{[\bi]} = 0$ for all~$i$ and~$j$.
Consider the~$6$-form
\be
\Lambda' = L\,\Omega 
          + 2\Re\left(L^{\bj\bk}_i\,\theta^i_\bj\w\omega_{[\bk]}\right).
\ee
The structure equations derived so far and the definitions made imply that%
\footnote{The reader familiar with the theory of Poincar\'e-Cartan
forms for first-order functionals will recognize this formula as
the first step in computing the so-called \emph{Betounes form} of the
Lagrangian~$\Lambda$.  While I could, at this point, compute the entire 
Betounes form, I will not need it in this article, so I leave this aside.  
For more information about the Betounes form, see~\cite{MR1985469}.}
\be
\d\Lambda' = 2\Re\left(\theta^i_\bj\w\Upsilon^\bj_i\right)
\ee
where
\be
\Upsilon^\bj_i 
= \d\bigl(L^{\bj\bk}_i\,\omega_{[\bk]}\bigr)
    +L^{\bj\bk}_\ell\,\phi^\ell_i\w\omega_{[\bk]}
    -L^{\bar\ell\bk}_i\,\overline{\phi^j_\ell}\w\omega_{[\bk]}\,.
\ee

Consequently, the Euler-Lagrange ideal~$\cE_\Lambda$ on~$J^1(\Js(M))$ 
for the first-order functional defined by~$\Lambda = L\,\Omega$ 
pulls back to~$\Fs^{(1)}$ to be generated by the $1$-forms~$\theta^i_\bj$,
their exterior derivatives, and the $6$-forms~$\Upsilon^\bj_i$ 
(whose real and imaginary parts constitute~$18$ independent forms).

\subsubsection{The Euler-Lagrange equations}
The upshot of these calculations is the following result, whose
proof has just been given in the course of the above discussion.

\begin{theorem}[Euler-Lagrange Equations]
\label{thm: EulerLagrange}
Let~$L(a) = c_1 |\det(a)|^2+c_2|P(\bar a,\bar a,{}^ta)|^2$
for~$a\in M_{3\times3}(\bbC)$ and real constants~$c_1$ and~$c_2$.
Define the polynomial functions~$L^{\bj\bk}_i(a)=-L^{\bk\bj}_i(a)$ 
as in~\eqref{eq: Ljkidefined}.  Then an almost complex structure~$J$
on~$M^6$ is a critical point of the functional~$F_c$ as defined
in~\eqref{eq: Fcdefined} if an only if, for any local $J$-linear
coframing~$\alpha:TU\to\bbC^3$ and $(0,2)$-connection~$\kappa$ satisfying
\be
\label{eq: 1ststreqalpha02}
\d\alpha= -\kappa\w\alpha + N(\alpha)\,
       \bpm \overline{\alpha^2\w\alpha^3}\strut\\
            \overline{\alpha^3\w\alpha^1}\strut\\
            \overline{\alpha^1\w\alpha^2}\strut \epm,
\ee
the following identities hold:
\be
\label{eq: ELeqns}
\d\bigl(L^{\bj\bk}_i(N(\alpha))\,\alpha_{[\bk]}\bigr)  
    +L^{\bj\bk}_\ell(N(\alpha))\,\kappa^\ell_i\w\alpha_{[\bk]}
    -L^{\bar\ell\bk}_i(N(\alpha))\,\overline{\kappa^j_\ell}\w\alpha_{[\bk]}
 = 0.
\ee
\end{theorem}

\begin{remark}[The nature of the equations]
As expected, Theorem~\ref{thm: EulerLagrange} shows that the
Euler-Lagrange equations for~$F_c$ are second-order and quasi-linear 
(i.e., linear in the second derivatives of~$J$).  The reader will also 
note that they constitute a system of $18$ equations.  
Of course, these equations can be expressed as the vanishing 
of a globally defined tensor on~$M$ 
and that will be done below in Corollary~\ref{cor: EulerLagrange},
but this is not particularly illuminating.  
The above formulation is good for practical calculations.

These equations cannot be elliptic since they are invariant 
under the diffeomorphism group.  It seems reasonable to expect 
the functional~$F_c$ to have `better' convexity properties 
when both~$c_1$ and~$c_2$ are positive, but I do not know 
whether this is the case.  However, an analysis that 
is somewhat long does show that, when~$c_1$ and~$c_2$ are positive, 
$F_c$ is elliptic for variations transverse to the action 
of the diffeomorphism group when linearized at any~$J$ 
for which~$\omega(J)$ is positive definite.  (Note that, when
$c_1$ and~$c_2$ are nonnegative, the functional~$F_c$
is bounded below by~$0$ and hence all integrable almost complex
structures are absolute minima for~$F_c$.  Of course, the 
integrable almost complex structures are critical points
of the functional~$F_c$ for any value of~$c$.)  

Finally, note that, even though~$\kappa$ is not uniquely determined
by~\eqref{eq: 1ststreqalpha02}, this ambiguity in~$\kappa$ does not
show up in the Euler-Lagrange equations~\eqref{eq: ELeqns}.  The
reason is that~\eqref{eq: 1ststreqalpha02} determines~$\kappa$ 
up to replacement by~${}^\ast\!\kappa$ where
\be
{}^\ast\!\kappa^i_j = \kappa^i_j + s^i_{jk}\,\alpha^k
\ee
for some functions~$s^i_{jk}=s^i_{kj}$.  Since~$\alpha_{[\bk]}$
contains~$\alpha^1\w\alpha^2\w\alpha^3$ as a factor,
${}^\ast\!\kappa^\ell_i\w\alpha_{[\bk]} = \kappa^\ell_i\w\alpha_{[\bk]}$,
so that the second term in the left hand side of~\eqref{eq: ELeqns}
is unaffected by the ambiguity in the choice of~$\kappa$.
It also follows that
\be
\overline{{}^\ast\!\kappa^i_\ell}\w\alpha_{[\bk]}
= \overline{\kappa^i_\ell}\w\alpha_{[\bk]} 
   + \overline{s^i_{\ell k}}\,\Omega
\ee
and this, coupled with the skewsymmetry~$L^{\bj\bk}_i=-L^{\bk\bj}_i$,
shows that the ambiguity in the choice of~$\kappa$ 
does not affect the final term of the left hand side of~\eqref{eq: ELeqns}.
\end{remark}

\begin{remark}[Degenerate critical points]
Of course, any~$J$ that satisfies~$\Phi(J)=0$ is a critical
point of the functional~$F_c$ when~$c_2=0$, just as any~$J$
that satisfies~$\psi(J)=0$ is a critical point of~$F_c$ when~$c_1=0$.

In particular, if~$d^{-1,2}:\cA^{1,0}\to\cA^{0,2}$
has rank at most~$1$ at each point, then one has both~$\Phi(J)=0$
and~$\psi(J)=0$, implying that all such almost complex structures
are critical for all of the~$F_c$.  

On the other hand, as will be seen below 
in Proposition~\ref{prop: supercritical}, 
the almost complex structure underlying a nearly K\"ahler structure 
will have~$\omega(J)>0$ and yet will be critical 
for all of the functionals~$F_c$.  
Thus, nondegenerate `supercritical' (i.e., critical for all~$F_c$) 
almost complex structures do occur.
\end{remark}

\subsubsection{The tensorial version}
Finally, I give the `tensorial' version of the Euler-Lagrange equations.
First, though, I note that because of the canonical isomorphism
\be
A^{p,q}(M) = A^{p,0}(M)\otimes A^{0,q}(M),
\ee
the $\pbar$-operators~$\pbar:A^{1,1}(M)\to A^{1,2}(M)$ and~$\pbar:A^{3,0}(M)
\to A^{3,1}(M)$ can be combined by the Leibnitz rule and exterior product
to yield a well-defined, first-order, linear operator
\be
\pbar: A^{1,1}(M)\otimes A^{3,0}(M)\to A^{1,0}(M)\otimes A^{3,2}(M).
\ee
Then Theorem~\ref{thm: EulerLagrange} can be expressed in the following
terms.

\begin{corollary}[The Euler-Lagrange equations]
\label{cor: EulerLagrange}
An almost complex structure~$J$ is critical for the functional~$F_c$ 
if and only if it satisfies
\be
\pbar\left(c_1\,E_1(J) 
+ c_2\bigl({\ts\frac13}E_2(J) + {\ts\frac23}E_0(J)\bigr)\right) = 0,
\ee
where the~$E_i(J)$ are defined as in~\eqref{eq: E1J}, \eqref{eq: E2J},
and \eqref{eq: E0J}.
\end{corollary}

\begin{proof}
This is an exercise in unraveling the definitions and applying 
Theorem~\ref{thm: EulerLagrange}.  The essential part of the calculation
is noting that the elementary identity
\be
\d P(a,b,c) = {\ts\frac13}\tr\bigl(Q(b,c)\,\d a+Q(c,a)\,\d b+Q(a,b)\,\d c\bigr)
\ee
can be expanded to yield a formula for~$\d L$ when
\be
L(a) = c_1\,\bigl|\det(a)\bigr|^2 + c_2\,\bigl|P(\bar a,\bar a,{}^ta)\bigr|^2.
\ee
Substituting this formula for~$\d L$ into the explicit equations of
Theorem~\ref{thm: EulerLagrange} and writing out the definition of~$\pbar$
in a coframing then yields the desired result.
\end{proof}

\begin{remark}[Reality and the tensors~$E_i(J)$]
\label{rem: reality&E_i}
By properties of the representation~$\rho$,
it is an invariant condition on~$J$ that, 
at each point, the Nijenhuis matrix~$N(\alpha)$ 
be a complex multiple of a Hermitian symmetric matrix.
  
In this case, which will be called 
the case of a \emph{Nijenhuis tensor of real type}, 
one clearly has
\be
\Phi(J) = {\ts\frac\iC8}\,\psi(J)\w\overline{\psi(J)}.
\ee
Moreover, by the formulae~\eqref{eq: E1J}, \eqref{eq: E2J}, 
and~\eqref{eq: E0J}, one also has
\be
E_1(J) = E_2(J) = E_0(J) = \omega(J)\otimes\psi(J).
\ee
In particular, by Corollary~\ref{cor: EulerLagrange}, 
an almost complex structure~$J$ with Nihenhuis tensor 
of real type is a critical point of all of the functionals~$F_c$ 
if and only if it is critical for some functional~$F_c$ 
where~$c = (c_1,c_2)$ satisfies~$c_1+c_2\not=0$.
In turn, this happens if and only if~$J$ satisfies the second-order condition
\be
\pbar\bigl(\omega(J)\,\otimes\,\psi(J)\bigl) = 0,
\ee
i.e.,
\be
\pbar\bigl(\omega(J)\bigl)\,\otimes\,\psi(J) 
+ \omega(J)\,\widehat\otimes\,\pbar\bigl(\psi(J)\bigr) = 0,
\ee
where the terms on the left hand side are regarded
as sections of~$A^{1,0}(M)\otimes A^{3,2}(M)$ via the canonical
isomorphism
\be
A^{1,2}(M)\otimes A^{3,0}(M) = A^{1,0}(M)\otimes A^{3,2}(M)
\ee
and the canonical homomorphism induced by wedge product
\be
A^{1,1}(M)\otimes A^{3,1}(M)\longrightarrow
A^{1,0}(M)\otimes A^{3,2}(M).
\ee
\end{remark}

\subsubsection{The nearly pseudo-K\"ahler case}
An interesting special case is the following one,
in which a second-order (determined) system of \textsc{pde} on~$J$
coupled with a first-order open condition on~$J$
implies the Euler-Lagrange equations for all of the functionals~$F_c$
and, moreover, a first-order closed condition on~$J$. 

\begin{proposition}[Nearly pseudo-K\"ahler structures]
\label{prop: nrlypsdoKahler}
Suppose that the almost complex structure~$J$ on~$M^6$ 
satisfies~$\pbar\omega=0$ while~$\psi$ is nonvanishing.
Then~$J$ is critical for all of the functionals~$F_c$ 
and its Nijenhuis tensor is of real type.
\end{proposition}

\begin{proof} Since~$\pbar\omega=0 = \overline{\p\omega}$, 
the relation~\eqref{eq: domega} simplifies to
\be
\d\omega = 3\Im(\psi) = 3{\ts\frac\iC2}(\overline{\psi}-\psi).
\ee
This implies~$\d\overline{\psi} = \d\psi$, which,
by type considerations, implies both that
\be
\pbar\psi = 0
\ee
and that~$d^{-1,2}\psi$ is a real-valued $(2,2)$-form.
Of course, this already implies that
\be
\pbar\bigl(\omega\otimes\psi) = 0.
\ee

Now, by \eqref{eq: dm12psi} and the formula for~$R(a)$, 
it follows from the reality of~$d^{-1,2}\psi$ that,
for any local complex coframing~$\alpha:TU\to\bbC^3$, the matrix
\be
R\bigl(N(\alpha)\bigr) = P\bigl(\overline{N(\alpha)},\overline{N(\alpha)},
{}^tN(\alpha)\bigr)\,N(\alpha),
\ee
is Hermitian symmetric.  Since
\be
\psi = P\bigl(\overline{N(\alpha)},\overline{N(\alpha)},
{}^tN(\alpha)\bigr)\,\alpha^1\w\alpha^2\w\alpha^3
\ee 
is assumed to be nonvanishing, it follows that~$N(\alpha)$
is a (nonzero) complex multiple of a Hermitian symmetric matrix,
i.e., that~$J$ has a Nijenhuis tensor of real type, as claimed.  

In particular, by Corollary~\ref{cor: EulerLagrange}, 
Remark~\ref{rem: reality&E_i}, 
and the formula~$\pbar\bigl(\omega\otimes\psi) = 0$, 
it now follows that~$J$ is critical for all~$F_c$, as claimed.
\end{proof}

\begin{remark}[Nomenclature]
In the case that~$\pbar\omega=0$ while~$\psi$ is nowhere vanishing,
one now sees by the formula for~$R(a)$ and the reality of
the Nijenhuis tensor that
\be
\d\psi = \d^{-1,2}\psi = 2\,\omega^2.
\ee
As will be seen in \S\ref{sec: ex&gen}, the equations~$\d\omega=3 \Im(\psi)$
and~$\d\psi = 2\,\omega^2$ characterize the so-called `strictly
nearly K\"ahler' case when~$\omega$ is a positive~$(1,1)$-form 
and~$\psi$ is nonvanishing.
However, as will be seen, there are examples of solutions of the above
equations for which~$\omega$ is not a positive~$(1,1)$-form, 
and this can be thought of as the `strictly nearly pseudo-K\"ahler' case.
\end{remark}

\section{Pseudo-holomorphic bundles}
\label{sec: pdoholbundles}

When~$(M,J)$ is a complex $n$-manifold, there is a close relationship
between holomorphic structures on complex bundles over~$M$ 
and connections on those bundles whose curvature satisfies 
certain restrictions.  This relationship is made somewhat
clearer in the context of Hermitian bundles over~$M$, so 
that is how the discussion will begin.

Suppose, first, that~$E\to M$ is a holomorphic bundle over~$M$.
Let~$\cH(E)$ denote the space of Hermitian metrics on~$E$.
As is well-known~\cite[p.~73]{MR0507725}, for any~$h\in\cH(E)$, 
there exists a unique connection~$\nabla^h$ on~$E$
that is compatible with~$h$ and satisfies~$(\nabla^h)^{0,1} 
= \pbar_E$.  The curvature~$(\nabla^h)^2$ of~$\nabla^h$
is an~$\End(E)$-valued $(1,1)$-form on~$M$.

Conversely, suppose that~$E$ is a complex vector bundle over~$M$,
that~$h\in\cH(E)$ is an Hermitian metric on~$E$, and that 
~$\nabla$ is an $h$-compatible connection on~$E$ whose
curvature~$\nabla^2$ is of type~$(1,1)$.  
Then there is a unique holomorphic structure on~$E$ such that
$\nabla^{0,1} = \pbar_E$.  

Of course, this correspondence between holomorphic structures
and connections with curvature of type~$(1,1)$ is not bijective,
precisely because there are many Hermitian structures on~$E$.
One way of tightening this correspondence is to impose some
extra conditions on the curvature form.  

For example, suppose that~$M$ carries a K\"ahler structure, 
represented by a K\"ahler form~$\eta$, say.%
\footnote{In what follows, one can allow the 
metric~$\eta$ to be indefinite (but always nondegenerate).
Thus, the~$(1,1)$-form~$\eta$ need not be positive.  
However, for the sake of brevity, I will usually only treat 
the definite case explicitly, only pointing out the places
where the sign makes a significant difference.}  
Then, for a given holomorphic bundle~$E$, 
an Hermitian metric~$h\in\cH(E)$
is said to be \emph{Hermitian-Yang-Mills} if it satisfies
\be
\tr_\eta\bigl((\nabla^h)^2\bigr) = \iC\lambda\,\id_E
\ee
for some constant~$\lambda\in\bbR$. 
This is a determined equation for~$h$ that is elliptic if~$\eta$
is definite and hyperbolic if~$\eta$ is indefinite. 

Alternatively, one can consider a complex bundle~$E\to M$
endowed with a fixed Hermitian metric~$h\in\cH(E)$ and look
for connections~$\nabla$ that are $h$-compatible and satisfy
\be
\bigl((\nabla)^2\bigr)^{0,2} = 0, \qquad
\tr_\eta\bigl((\nabla)^2\bigr)^{1,1} = \iC\lambda\,\id_E
\ee
for some constant~$\lambda\in\bbR$.  

For the importance of these concepts in algebraic and complex differential
geometry, see~\cite{MR0765366} and~\cite{MR0861491}.
This introductory discussion is meant to motivate the following definitions.

\begin{definition}[Pseudo-holomorphic bundles]
When~$(M,J)$ is an almost complex manifold and $E\to M$ is
a complex vector bundle over~$M$ endowed with a Hermitian
metric~$h$, an $h$-compatible connection~$\nabla$ on~$E$
is said to define a \emph{pseudo-holomorphic} structure on~$E$
if~$\nabla^2$ has curvature of type~$(1,1)$.  

If, in addition,~$\eta$ is a (possibly indefinite) 
almost Hermitian metric on~$(M,J)$,
and~$\tr_\eta(\nabla^2) = \iC\lambda\,\id_E$, the connection~$\nabla$
is said to be \emph{pseudo-Hermitian-Yang-Mills}.
\end{definition}

\begin{remark}[Philosophy]
The main reason that this generalization of holomorphic bundle
and Hermitian-Yang-Mills connections might be of interest is
that these can exist and be nontrivial even when~$(M,J)$ is
not an integrable complex manifold.  In fact, the main
point of this section is to point out that there is a rather large class
of almost complex structures in (real) dimension~$6$ that, 
at least locally, support `as many' pseudo-Hermitian-Yang-Mills 
connections as the integrable complex structures do.

These pseudo-Hermitian-Yang-Mills connections on almost
complex manifolds, particularly the $\Gtwo$-invariant $6$-sphere, 
also show up in the study of singular behaviour of Yang-Mills
connections.  See~\S5.3 of \cite{Tian2000} for a further
discussion of this.

Also, while it is not, strictly speaking, necessary to restrict 
attention to connections~$\nabla$ compatible with a Hermitian
structure on~$E$, it is convenient and it simplifies the theory
in some respects.  I will leave it to the reader to determine
when this restriction can be dropped.  In any case, since~$E$
is only being regarded as a complex bundle (rather than a 
holomorphic one), any two Hermitian structures on~$E$ are equivalent 
under the action of the smooth gauge group~$\Aut(E)$, so the 
specification of a Hermitian norm on~$E$ is not a delicate issue.
\end{remark}

\begin{proposition}[Curvature restrictions]
\label{prop: curvrestr}
If~$\nabla$ is a pseudo-holomorphic structure on~$E\to M$,
then~$\nabla^2 = (\nabla^2)^{1,1}$ takes values in~$\Hom(E)\otimes\cK$
where
\be
\cK=\ker\left(\d^{-1,2}+\d^{2,-1}:\cA^{1,1}\to\cA^{0,3}\oplus\cA^{3,0}\right).
\ee
\end{proposition}

\begin{proof}
Let~$r$ be the rank of~$E$ and
let~$\us = (\us_1,\ldots,\us_r)$ be a unitary basis of~$\Gams(U,E)$
for some open set~$U\subset M$.  Then $\nabla\us = \us A$,
where~$A$ is a $1$-form on~$U$ with values in~$\euu(r)$.
By hypothesis, the curvature form~$F = \d A + A\w A$ is a $(1,1)$-form
with values in~$\euu(r)$.  By the Bianchi identity
\be
\d F = F\w A - A \w F.
\ee
Since~$F$ is a $(1,1)$-form, the $\cA^{0,3}\oplus\cA^{3,0}$-component 
of the right hand side of this equation vanishes.  Thus,
\be
\bigl(\d^{-1,2}+\d^{2,-1}\bigr)(F) = 0,
\ee
as desired.
\end{proof}

\begin{remark}[Generic triviality]
When~$n$ is sufficiently large and~$J$ (or more, specifically, the
Nijenhuis tensor of~$J$) is sufficiently generic, the map~$d^{-1,2}+d^{2,-1}$
is injective on~$\cA^{1,1}$.  In such a case, any pseudo-holomorphic
structure~$\nabla$ on~$E$ is necessarily a flat connection, so that~$E$
is a flat bundle.  This is analogous to the fact that the sheaf of
holomorphic functions on a general nonintegrable almost complex structure 
is just the constant sheaf.  
\end{remark}

The interesting case is going to be when there are plenty of non-flat 
pseudo-holomorphic structures.  Note that the set of complex bundles
over~$M$ that admit pseudo-holomorphic structures is closed under 
sums, products, and taking duals.  It would appear that a good way 
to find almost complex structures that support many pseudo-holomorphic
bundles would be to look for ones that support many pseudo-holomorphic
line bundles.  It is to this case that I will now turn.

\subsection{The line bundle case}
From now on, I will restrict attention to the $6$-dimensional case.

Consider the case of a line bundle~$E\to M$.  
Then any pseudo-holomorphic structure~$\nabla$ on~$E$
has curvature of the form~$\nabla^2 = \iC\phi\otimes\id_E$
where~$\phi$ is a closed, real-valued~$(1,1)$-form. 

Using a local complex coframing~$\alpha = (\alpha^i)$ on
the open set~$U\subset M$, write
\be
\phi = {\ts\frac{\iC}2} F_{i\bj}\,\alpha^i\w\overline{\alpha^j}
\ee
where~$F_{i\bj}=\overline{F_{j\bi}}$ can be regarded as a Hermitian
symmetric $3$-by-$3$ matrix~$F$ defined on~$U$.  One computes
\be
\label{eq: dm12phi}
\d^{-1,2}\phi = {\ts\frac{\iC}2}\tr\bigl({}^tN(\alpha)F\bigr)\,
\overline{\alpha^1\w\alpha^2\w\alpha^3}.
\ee
Thus, closure implies two linear zeroth-order equations on~$\phi$, 
namely, the real and imaginary parts 
of the equation~$\tr\bigl({}^tN(\alpha)F\bigr)=0$. 

This algebraic system combined with the already overdetermined
system~$\d\phi=0$ for~$\phi\in\cA^{1,1}$ is generally
not involutive.%
\footnote{In fact, by a somewhat long argument that
involves higher order invariants of the almost complex 
structure, it can be shown that, for a sufficiently generic
almost complex structure~$J$ on a $6$-manifold~$M$, the sheaf of
closed $(1,1)$-forms on~$M$ is trivial.  Since the methods
of this argument will not be used further in this note, I 
will not give details here.}
However, there is a class of almost complex structures
characterized by a first-order condition that generalizes
the integrable case and for which the system~$\d\phi=0$
for~$\phi\in\cA^{1,1}$ is involutive.  It is to this
condition that I now turn.

Note that, because~$F$ is Hermitian-symmetric, 
the equation~$\tr\bigl({}^tN(\alpha)F\bigr)=0$ will reduce 
to a single linear equation on~$F$ if and only if~$N(\alpha)$ 
is a complex multiple of a Hermitian-symmetric matrix, i.e.,
if and only if~$J$ has a Nijenhuis tensor of real type
(see Remark~\ref{rem: reality&E_i}).

It turns out that simply having its Nijenhuis tensor be
of real type is a rather awkward condition on~$J$
because the union of the $\rho$-orbits 
of the Hermitian symmetric matrices 
in~$M_{3\times3}(\bbC)$ has complicated singularities.
It is better to work with a somewhat stronger condition
that `tames' this singularity, which motivates the following 
definition:  

\begin{definition}[Quasi-integrability]
An almost complex structure~$J$ on~$M^6$
is said to be \emph{elliptically quasi-integrable} if, for any 
local $J$-linear coframing~$\alpha$, the matrix~$N(\alpha)$ takes
values in the $\rho$-orbit of~$I_3\in M_{3\times3}(\bbC)$
union with the zero element.  It is said to be \emph{strictly
elliptically quasi-integrable} if, for any local $J$-linear 
coframing~$\alpha$, the matrix~$N(\alpha)$ takes values 
in the $\rho$-orbit of~$I_3\in M_{3\times3}(\bbC)$.

Similarly, an almost complex structure~$J$ on~$M^6$
is said to be \emph{hyperbolically quasi-integrable} if, for any 
local $J$-linear coframing~$\alpha$, the matrix~$N(\alpha)$ takes
values in the $\rho$-orbit of~$\diag(1,-1,-1)\in M_{3\times3}(\bbC)$
union with the zero element.  It is said to be \emph{strictly
hyperbolically quasi-integrable} if, for any local $J$-linear 
coframing~$\alpha$, the matrix~$N(\alpha)$ takes values 
in the $\rho$-orbit of~$\diag(1,-1,-1)\in M_{3\times3}(\bbC)$.
\end{definition}

\begin{remark}[Nonstrictness]
The main reason for including the possibility of the Nijenhuis
tensor vanishing in the definition of quasi-integrability
is so that integrable almost complex structures will be 
quasi-integrable.  The main focus of the rest of this section
will be the strictly quasi-integrable case.

The justification for the term `quasi-integrable' will appear below.  
Basically, the reason is that the overdetermined system 
describing pseudo-Hermitian-Yang-Mills connections,
which would not be involutive 
for the general almost complex structure in dimension~$6$, 
turns out to be very well-behaved in the quasi-integrable case.
\end{remark}

\begin{proposition}[Structure reduction]
If~$J$ is an elliptic strictly quasi-integrable almost complex structure 
on the $6$-manifold~$M$, then there is a canonical $\SU(3)$-structure
on~$M$ whose local sections consist of the local $J$-linear coframings
$\alpha:TU\to\bbC^3$ for which~$N(\alpha) = I_3$.
The invariant~$(1,1)$-form~$\omega = \omega(J)$
is a positive~$(1,1)$-form on~$M$ while $\psi$ is a nonvanishing
$(3,0)$-form on~$M$.

Similarly, if~$J$ is a hyperbolic strictly quasi-integrable 
almost complex structure on the $6$-manifold~$M$, then there is 
a canonical $\SU(1,2)$-structure
on~$M$ whose local sections consist of the local $J$-linear coframings
$\alpha:TU\to\bbC^3$ for which~$N(\alpha) = \diag(1,-1,-1)$.
The invariant~$(1,1)$-form~$\omega = \omega(J)$
is a~$(1,1)$-form of signature~$(1,2)$ on~$M$ while 
$\psi$ is a nonvanishing $(3,0)$-form on~$M$.

In either case \oppar hyperbolic or elliptic\clpar, these forms satisfy
\be
{\ts\frac16}\,\omega^3 = {\ts\frac\iC8}\,\psi\w\overline\psi = \Phi>0.
\ee
\end{proposition}

\begin{proof}
The elliptic case follows from the evident fact
that the $\rho$-stabilizer of~$I_3\in M_{3\times3}(\bbC)$ 
is~$\SU(3)$ while~$P(I_3,I_3,I_3) = 1$ and~$Q(I_3,I_3) = I_3$.
The hyperbolic case is similar since, if~$D = \diag(1,-1,-1)$,
then $P(D,D,D) = 1$ and~$Q(D,D) = D$.
\end{proof}

Recall that a $\U(3)$-structure 
(also called an \emph{almost Hermitian structure}) on~$M^6$ 
is a pair~$(J,\eta)$ consisting of an almost complex
structure~$J$ on~$M$ and a $2$-form~$\eta$ that is a positive $(1,1)$-form
with respect to~$J$.  A $J$-linear coframing~$\alpha:TU\to\bbC^3$
defined on~$U\subset M$ is said to be~\emph{$(J,\eta)$-unitary}
if~$\alpha$ defines a unitary isomorphism between~$T_xU$ and~$\bbC^3$
(with its standard Hermitian structure) for all~$x\in U$.
One has a similar description of~$\U(1,2)$-structures on~$M$.
The transition rule~\eqref{eq: Ntransition} then implies that
the following definitions are meaningful.

\begin{definition}[Quasi-integrable unitary structures]
A $\U(3)$-structure~$(J,\eta)$ on~$M^6$ is said to be \emph{quasi-integrable} 
if any local~$(J,\eta)$-unitary coframing $\alpha:TU\to\bbC^3$ 
satisfies~$N(\alpha) = \lambda I_3$ for some function~$\lambda:U\to\bbC$.

Similarly, a~$\U(1,2)$-structure~$(J,\eta)$ on~$M^6$ is said
to be \emph{quasi-integrable} if any local~$(J,\eta)$-unitary coframing 
$\alpha:TU\to\bbC^3$ satisfies~$N(\alpha) 
= \diag(\lambda,-\lambda,-\lambda)$ for some function~$\lambda:U\to\bbC$.
\end{definition}  

\begin{remark}[Relations]
If~$(J,\eta)$ is quasi-integrable on~$M$, then~$J$ is quasi-integrable
on~$M$.  However, there do exist quasi-integrable~$J$ for which
there does not exist a corresponding~$\eta$ such that~$(J,\eta)$
is quasi-integrable.  

The difficulty is caused by the places where
the Nijenhuis tensor of~$J$ vanishes.  If~$J$ is quasi-integrable,
then on the open set~$M^*\subseteq M$ where its Nijenhuis tensor 
is nonzero (i.e., where~$J$ is strictly quasi-integrable), 
the pair~$\bigl(J,\omega(J)\bigr)$ is quasi-integrable. In fact,
any~$\eta$ on~$M^*$ such that~$(J,\eta)$ is quasi-integrable 
is necessarily a multiple of~$\omega(J)$.  However, it can
happen that one cannot smoothly extend any such positive~$\eta$
as a positive $(1,1)$-form on all of~$M$.  See Remark~\ref{rem: sing}.
\end{remark}

\begin{remark}[The quasi-integrable condition]
Because the orbit~$\rho\bigl(\GL(3,\bbC)\bigr){\cdot}I_3\subset 
M_{3\times3}(\bbC)$ is diffeomorphic to~$\GL(3,\bbC)/\SU(3)$,
this orbit has dimension~$10$ and codimension~$8$
in~$M_{3\times3}(\bbC)$.  Thus, the condition of strict
elliptic quasi-integrability
can be thought of as a system of $8$ first-order equations (plus
some inequalities due to the fact that the orbit is not closed)
for the almost complex structure~$J$.  The case of strict hyperbolic
quasi-integrability is similar. The nature of these
equations will be discussed below.

Somewhat better behaved is the quasi-integrability condition 
for a pair~$(J,\eta)$.  In this case, quasi-integrability is 
a closed condition and is easily seen to be a system of~$16$
first-order equations for the $\U(3)$- or ~$\U(1,2)$-structure~$(J,\eta)$.
Since quasi-integrability is unaffected by 
replacing~$(J,\eta)$ by~$(J,\mu\,\eta)$ for any smooth
positive~$\mu$, it would be more natural to formulate 
this concept for $\bbR^+{\cdot}\U(3{-}q,q)$-structures, 
rather than~$\U(3{-}q,q)$-structures, 
but I will wait until the discussion 
in~\S\ref{ssec: generality} to do this.
\end{remark}

\begin{remark}[Other nomenclatures]
\label{rem: GrayHervilla}
It should also be noted that quasi-integrability 
for a $\U(3)$-structure is what Gray and Hervella~\cite{MR0581924}
refer to as type~$\cW_1{\oplus}\cW_3{\oplus}\cW_4$
(or, alternatively~$\cG_1$, see~\cite{MR0431008,MR0577582}), 
since, quasi-integrability for a~$\U(3)$-structure 
is the vanishing of the $\U(3)$-torsion component 
that Gray and Hervella denote by~$W_2$.  
\end{remark}

\begin{proposition}
\label{prop: qiimpliespdoHYM}
Let~$(J,\eta)$ be a strictly quasi-integrable $\U(3{-}q,q)$-structure 
on~$M^6$.  Then any pseudo-holomorphic stucture~$\nabla$ on an Hermitian 
bundle~$E\to M$ is pseudo-Hermitian-Yang-Mills.  In fact,
\be
\tr_\eta\bigl(\nabla^2\bigr) = 0.
\ee
\end{proposition}

\begin{proof}
Under these conditions, the formula~\eqref{eq: dm12phi}
implies that the kernel space~$\cK$ 
defined in Proposition~\ref{prop: curvrestr} 
is equal to the set of~$\phi\in\cA^{1,1}$ such that~$\tr_\eta(\phi)=0$.
\end{proof}

The main interest in quasi-integrability comes from the following
elliptic result and its hyperbolic analog, which ensure that, at least locally, 
there are many pseudo-holomorphic line bundles on quasi-integrable manifolds.

\begin{proposition}[The fundamental exact sequence]
\label{prop: fundextseq}
Let~$(J,\eta)$ be a quasi-integrable $\U(3)$-structure on~$M^6$
and let
\be
\cA^{1,1}_0(M) = \left\{\phi\in \cA^{1,1}(M)\mid \tr_\eta(\phi)=0\right\},
\ee
be the $\eta$-primitive part of the $(1,1)$-forms on~$M$.
Let~$\cZ^{1,1}_0(M)\subset \cA^{1,1}_0(M)$ denote the closed
forms in~$\cA^{1,1}_0(M)$.  

Then~$d^{-1,2}+d^{2,-1}$ vanishes on~$\cA^{1,1}_0$, the subcomplex
\be
\label{eq: fundsubcomplex}
0\longrightarrow \cA^{1,1}_0\longrightarrow \cA^{2,1}{\oplus}\cA^{1,2}
\longrightarrow\cA^4\longrightarrow\cA^5\longrightarrow\cA^6\longrightarrow0
\ee
of the deRham complex on~$M$ is a locally exact elliptic complex,
and, consequently, the sheaf complex
\be
\label{eq: fundres}
0\longrightarrow\cZ^{1,1}_0\longrightarrow \cA^{1,1}_0
\longrightarrow \cA^{2,1}{\oplus}\cA^{1,2}
\longrightarrow\cA^4\longrightarrow\cA^5\longrightarrow\cA^6\longrightarrow0
\ee
is a fine resolution of~$\cZ^{1,1}_0$.
\end{proposition}

\begin{proof}
To begin, note that~$\cA^{1,1}(M)$ is the space of sections of 
a (complex) vector bundle of rank~$9$ over~$M$ and that $\cA^{1,1}_0(M)$ 
is the space of sections of a (complex) subbundle of rank~$8$.
If~$U\subseteq M$ is an open subset on which there exists
a $(J,\eta)$-unitary coframing~$\alpha:TU\to\bbC^3$, then one
has~$N(\alpha)= \lambda I_3$ for some~$\lambda\in C^\infty(U)$
and, moreover, for any~$\phi\in\cA^{1,1}_0(U)$, one has
\be
\phi ={\ts\frac\iC2}\, F_{i\bj}\,\alpha^i\w\overline{\alpha^j}
\ee
for some $3$-by-$3$ matrix of functions~$F = \bigl(F_{i\bj}\bigr)$
satisfying~$\tr(F)=0$.  Of course, by the formula~\eqref{eq: dm12phi},
it follows that~$(\d^{-1,2}+\d^{2,-1})(\phi)=0$.  

Thus,
$\d\bigl(\cA^{1,1}_0(M)\bigr)$ lies in~$\cA^{2,1}(M)\oplus\cA^{1,2}(M)$
and it follows that \eqref{eq: fundsubcomplex} is indeed a subcomplex
of the (complexified) deRham complex.  (Note, that, moreover, each
of the vector spaces involved is the complexification of the real
subspace consisting of the real-valued forms in that subspace.)

Since it is a subcomplex of the deRham complex,
for any nonzero (real) cotangent vector~$\xi$, 
the symbol~$\sigma_\xi$ of the operators 
in the sequence~\eqref{eq: fundsubcomplex} 
is simply left exterior multiplication by~$\xi$.  
This symbol sequence is obviously exact at~$\cA^5$
and~$\cA^6$. Moreover, elementary linear algebra shows that~$\cA^{1,1}_0$
contains no nonzero decomposable $2$-forms, so the symbol $\sigma_\xi$
is injective at~$\cA^{1,1}_0$.  By rank count (the bundles
have ranks~$8$, $18$, $15$, $6$ and~$1$ in order), then, 
it follows that the symbol sequence will be exact 
at~$\cA^{2,1}(M)\oplus\cA^{1,2}(M)$ if and only if it is exact
at~$\cA^4(M)$.  To prove that it is exact at~$\cA^4$, it suffices
to show that any~$4$-form of the form~$\xi\w\Upsilon$ can be 
written in this form where~$\Upsilon$ 
lies in~$\cA^{2,1}(M)\oplus\cA^{1,2}(M)$.  However, this is immediate
from exterior algebra considerations.  
Thus, the sequence~\eqref{eq: fundsubcomplex} is an elliptic complex.

To prove local exactness is somewhat more subtle and I will not give 
full details here.  For the results needed in the following argument and
the definitions of the terms involved
see~\cite{MR1083148}, especially Chapter~10, \S3.
(Note, however, that local exactness is obvious at the
places~$\cA^5$ and~$\cA^6$ since this is just the Poincar\'e Lemma.
The issue is local exactness at the other places.)

First of all, a calculation shows that
the Cartan characters of the overdetermined system~$\d\phi=0$
for~$\phi\in\cA^{1,1}_0$ are
\be
\label{eq: fundsubcomplexCartanchars}
(s_0,s_1,s_2,s_3,s_4,s_5,s_6) = (0,0,2,4,6,4,0)
\ee
and that this system is indeed involutive.  It then follows
that the Spencer resolution of the overdetermined operator
$\d:\cA^{1,1}_0\to\cA^{2,1}\oplus\cA^{1,2}$ is formally exact
and a computation verifies that~\eqref{eq: fundsubcomplex}
is, in fact, this Spencer resolution.  

Finally, the  operator~$\d:\cA^{1,1}_0\to\cA^{2,1}\oplus\cA^{1,2}$ 
is seen to satisfy Singer's $\delta$-estimate, which, together
with the formal exactness already derived, implies 
that~\eqref{eq: fundsubcomplex} is indeed locally exact
(see \cite[Theorem 3.10]{MR1083148}).
\end{proof}

\begin{remark}[The hyperbolic case]
\label{rem: hyperbolicfundlemma}
In the case of a hyperbolic quasi-integrable~$(J,\eta)$, 
the complex~\eqref{eq: fundsubcomplex} is still formally exact locally.
After all, the overdetermined system~$\d\phi=0$ for~$\phi\in\cA^{1,1}_0$
is still involutive and has Cartan characters 
given by~\eqref{eq: fundsubcomplexCartanchars}.
However, the local exactness is not obvious 
since one does not have an analog of the $\delta$-estimate.  
\end{remark}

\begin{remark}[Relation with Hodge theory]
\label{rem: Hodge}
Note that, when~$M^*\subseteq M$, the open set on which~$J$
is strictly quasi-integrable, is dense in~$M$, the space~$\cA^{1,1}_0(M)$
is equal to the kernel of~$\d^{-1,2}+\d^{2,-1}$ on~$\cA^{1,1}(M)$.
In particular, in this case, the vector space~$\cZ^{1,1}_0(M)$ 
is equal to the space of closed $(1,1)$-forms on~$M$.

By Proposition~\ref{prop: fundextseq}, when $M$ is compact
and~$(J,\eta)$ is elliptically quasi-integrable,
the space~$\cZ^{1,1}_0(M)$ is finite dimensional 
and consists of the closed `primitive' $(1,1)$-forms
on~$M$.  In general, one does not know that these forms are
harmonic in the usual sense.  The reason is that, 
when $\phi\in\bigl(\cA^{1,1}_0\bigr)_\bbR$, one has
\be
\ast\phi = -\phi\w\eta,
\ee
where~$\ast$ is the Hodge star operator for the underlying metric
and orientation associated to the~$\U(3)$-structure~$(J,\eta)$. 
When~$\phi$ is closed, this only gives
\be
\d(\ast\phi) = -\phi\w\d\eta = \phi\w(\p\eta + \pbar\eta),
\ee
so that, unless~$(J,\eta)$ satisfies~$\pbar\eta=0$ (which is not 
automatic), one will not generally have $\phi$ be coclosed.

In any case, if~$\phi=\bar\phi\in\cZ^{1,1}_0(M)$ satisfies the
condition that~$[\phi]\in H^2(M,\bbR)$ is an integral cohomology
class, then there will exist a complex line bundle~$E$ endowed
with an Hermitian structure~$h$ and an~$h$-compatible 
connection~$\nabla$ such that~$\nabla^2 = 2\pi\iC\phi\otimes\id_E$.
In particular, $\nabla$ will define a pseudo-holomorphic
structure on~$E$, one that is, moreover, pseudo-Hermitian-Yang-Mills.
\end{remark}

\begin{remark}[Potentials]
As the reader will recall, in the integrable case, a closed, 
smooth $(1,1)$-form~$\phi$ can be locally written in the 
form~$\phi = \iC\,\p\pbar f$ for some smooth `potential' function~$f$.  
Moreover, when~$(J,\eta)$ is Hermitian (i.e., $J$ is integrable),
an $\eta$-primitive~$(1,1)$ form is locally expressible
as~$\iC\,\p\pbar f$ where~$f$ is harmonic with respect to the 
underlying metric.  

Unfortunately, in the non-integrable case, there does not
appear to be such a local potential formula for the closed 
$(1,1)$-forms.
\end{remark}

\subsection{General rank}
In fact, in the quasi-integrable case, the above results for
line bundles generalize in a natural way to all ranks:

\begin{proposition}
\label{prop: ellipticqipdoHYMstructs}
Let~$(J,\eta)$ be a quasi-integrable structure on~$M^6$
and let~$E\to M$ be a rank~$r$ complex bundle endowed with 
an Hermitian structure~$h$.  The overdetermined system
\be
\label{eq: pdoHYMeqs}
(\nabla^2)^{0,2}=\tr_\eta(\nabla^2) = 0
\ee
for $h$-compatible connections~$\nabla$ on~$E$ is involutive,
with Cartan characters
\be
\label{eq: pdoHYMeqschars}
(s_0,s_1,s_2,s_3,s_4,s_5,s_6) = (0,0,0,0,2r^2,3r^2,r^2).
\ee
\end{proposition}

\begin{proof}
I will give the proof in the case that~$(J,\eta)$ is a $\U(3)$-structure.
The $\U(1,2)$-structure case is entirely analogous.

The claimed result is local, so let~$U\subset M$ be an open set over which
both~$U$ and~$T^{1,0}(M)$ are trivial as complex bundles.  Let~$\us
= (\us_1,\ldots,\us_r)$ be an $h$-unitary basis of the section 
of~$E$ over~$U$ and let~$\alpha:TU\to\bbC^3$ be a $(J,\eta)$-unitary
coframing on~$U$.  Then there exist $1$-forms~$\phi^i_j$ and
a smooth function~$\lambda$ on~$U$ such that~$\alpha$ satisfies
\be
\d\alpha^i = -\phi^i_l\w\alpha^l + \lambda\,\overline{\alpha^j\w\alpha^k}
\ee
where~$(i,j,k)$ is any even permutation of~$(1,2,3)$.  In particular,
note that
\be
\label{eq: dalphaformula}
\d(\alpha^1\w\alpha^2\w\alpha^3) 
= -\tr(\phi)\w \alpha^1\w\alpha^2\w\alpha^3
    + 2\lambda\,\eta^2,
\ee
since~$\eta = {\ts\frac\iC2}\,\bigl(\alpha^1\w\overline{\alpha^1}
+ \alpha^2\w\overline{\alpha^2} + \alpha^3\w\overline{\alpha^3}\bigr)$.

Now, any $h$-compatible connection~$\nabla$ on~$E$ over~$U$ 
will be of the form
\be
\nabla\us = \us \otimes \bigl(\pi-{}^t\bar\pi\bigr)
\ee
where
\be
\pi = p_1\,\alpha^1 + p_2\,\alpha^2 + p_3\,\alpha^3 
\ee
for some~$p_i:U\to M_{r\times r}(\bbC)$.   Letting
\be
\label{eq: Fdefined}
F = \d(\pi-{}^t\bar\pi\bigr) + (\pi-{}^t\bar\pi\bigr)\w(\pi-{}^t\bar\pi\bigr)
= -{}^t\overline{F}
\ee
denote the curvature of~$\nabla$ as usual, the condition 
that~$F = F^{1,1}$ can be written as
\be
F\w\alpha^1\w\alpha^2\w\alpha^3 = 0
\ee
and the condition that~$\tr_\eta(F) = 0$ can be written in the form
\be
F\w\eta^2 = 0.
\ee

This motivates the following construction:
Let~$X = U\times \bigl(M_{r\times r}(\bbC)\bigr)^3$, and regard
the projections~$p_i:X\to M_{r\times r}(\bbC)$ as vector-valued
functions on~$X$.  Define~$\pi = p_1\,\alpha^1 + p_2\,\alpha^2 + p_3\,\alpha^3$
and regard~\eqref{eq: Fdefined} as defining a $\euu(r)$-valued
$2$-form~$F$ on~$X$.  Let~$\cI$ be the ideal on~$X$ algebraically generated by 
the~$2r^2$ components 
of the $5$-form~$\Upsilon = F\w\alpha^1\w\alpha^2\w\alpha^3$
and the~$r^2$ components of the $6$-form~$\Psi = F\w\eta^2$.   
Note that, by the Bianchi identity and~\eqref{eq: dalphaformula},
one has the crucial identity
\be
\begin{aligned}
\d\Upsilon &= \d F \w \alpha^1\w\alpha^2\w\alpha^3 
               + F\w d\bigl(\alpha^1\w\alpha^2\w\alpha^3 \bigr)\\
  &= \bigl(F\w (\pi-{}^t\bar\pi) - (\pi-{}^t\bar\pi)\w F\bigr) \w 
       \alpha^1\w\alpha^2\w\alpha^3 \\
  &\qquad\qquad      + F\w \bigl(-\tr(\phi)\w \alpha^1\w\alpha^2\w\alpha^3
                   + 2\lambda\,\eta^2\bigr)\\
  &= -\Upsilon\w(\pi-{}^t\bar\pi) - (\pi-{}^t\bar\pi)\w\Upsilon
       -\tr(\phi)\w\Upsilon + 2\lambda\,\Psi,
\end{aligned}
\ee 
which shows that~$\cI$ is differentially closed.  (It is for this
that the hypothesis of quasi-integrability is needed.)

The integral manifolds of~$\bigl(\cI,\,\iC\,\alpha^1\w\alpha^2\w\alpha^3\w
\overline{\alpha^1\w\alpha^2\w\alpha^3}\bigr)$ in~$X$ are, by construction,
locally the graphs of the $h$-compatible connections on~$E$ over~$U$
that satisfy~\eqref{eq: pdoHYMeqs}.  The verification that this ideal
with independence condition is involutive and has characters as described 
in~\eqref{eq: pdoHYMeqschars} is now routine.
\end{proof}

\begin{remark}[Gauge fixing]
The reader should not be surprised that the last nonzero character
in~\eqref{eq: pdoHYMeqschars} is $s_6 = r^2$ rather than~$0$.  This
is to be expected since the equations~\eqref{eq: pdoHYMeqs} are
gauge invariant.  In terms of the local representation defined above,
this implies that if~$\pi = p_i\,\alpha^i$ does satisfy these equations,
then, for any smooth map~$g:U\to\U(r)$, 
the $\eugl(r,\bbC)$-valued $(1,0)$-form
\be
\pi^g = g^{-1}\,\p g + g^{-1}\pi g
\ee
will also satisfy these equations.  As usual, one can break this
gauge-invariance by imposing an additional equation, such as
the `Coulomb gauge' condition
\be
\d(\pi+{}^t\bar\pi)\w\eta^2 = 0,
\ee
and the reader can check that this yields a system whose characters
are
\be
\label{eq: pdoHYMeqscharsgauged}
(s_0,s_1,s_2,s_3,s_4,s_5,s_6) = (0,0,0,0,2r^2,4r^2,0).
\ee
Moreover, when~$\eta$ is positive definite, 
this is an elliptic system, which implies, in particular,
that one has elliptic regularity for Coulomb-gauged 
pseudo-Hermitian-Yang-Mills connections 
over a quasi-integrable~$\U(3)$-structure. 

In particular, note that this shows that pseudo-Hermitian-Yang-Mills 
connections have the same degree of local generality for
quasi-integrable~$\U(3)$-structures as for integrable~$\U(3)$-structures.
In the case of real-analytic quasi-integrable~$\U(3)$-structures, 
the Coulomb-gauged pseudo-Hermitian-Yang-Mills connections are
also real-analytic.

It seems very likely, though I have not tried to verify this, 
that there would be analogs of the Uhlenbeck removable singularities 
theorem generalizing to the quasi-integrable case 
the ones known in the integrable case~\cite{MR0765366,Tian2000,MR0861491}.

Naturally, in the hyperbolic case, i.e., for quasi-integrable
$\U(1,2)$-structures, one does not have elliptic regularity.
\end{remark}

\section{Examples and Generality}
\label{sec: ex&gen}

In this section, after some examples have been given of
(strictly) quasi-integrable almost complex structures, 
the question of how general such structures are will be addressed.

\subsection{The $6$-sphere}
\label{ssec: 6sphere}
The most familiar example of of a non-integrable almost complex
structure is the $\Gtwo$-invariant almost complex structure on
the~$6$-sphere.

The group~$\Gtwo$ is the subgroup of~$\GL(7,\bbR)$
that preserves the $3$-form
\be
\phi
=\d x^{123}+\d x^{145}+\d x^{167}+\d x^{246}-\d x^{257}-\d x^{356}-\d x^{347},
\ee
where~$x = (x^i)$ are the standard coordinates on~$\bbR^7$ 
and~$\d x^{ijk} = \d x^i\w\d x^j\w\d x^k$.  It is not difficult to show
that~$\Gtwo$ is a subgroup of~$\SO(7)$, i.e., that it preserves the
standard inner product and orientation on~$\bbR^7$.  It is connected
and has dimension~$14$.  It acts transitively on~$S^6\subset\bbR^7$
and the stabilizer of a point in~$S^6$ is isomorphic to~$\SU(3)$, 
i.e.,~$\Gtwo\cap\SO(6) = \SU(3)$.  

In particular, it follows that~$\Gtwo$ preserves a $\SU(3)$-structure 
on~$S^6$.  One can define this structure by defining its associated
invariant differential forms.  The invariant $(1,1)$-form~$\omega$
is defined to be the pullback to~$S^6$ of the $\Gtwo$-invariant 
$2$-form~$R\lhk\phi$, where~$R$ is the radial vector field on~$\bbR^7$. 
The invariant~$(3,0)$-form~$\psi$ is defined to be the pullback to~$S^6$
of the complex-valued $3$-form~$R{\lhk}(\ast\phi)+\iC\phi$, 
where~$\ast\phi$ is the $4$-form that is Hodge dual 
to~$\phi$ in~$\bbR^7$. The Euler identities
\be
\d(R\lhk\phi) = 3\,\phi\qquad\text{and}\qquad 
\d\bigl(R\lhk(\ast\phi)\bigr) = 4\,{\ast}\phi
\ee
together with a little algebra then imply the relations
\be
\d\omega = 3\Im(\psi)  \qquad\text{and}\qquad 
\d\psi = 2\,\omega^2.
\ee
(The almost complex structure~$J$ is defined to be the unique
one for which~$\psi$ is of type~$(3,0)$.)

Let~$u:\Gtwo\to S^6$ be defined by~$u(g) = g{\cdot}u_0$ for some
fixed unit vector~$u_0\in S^6$. It is then easy to verify that 
on~$\Gtwo$ there exist complex-valued left-invariant 
forms~$\alpha_i$ and~$\kappa_{i\bj} = -\overline{\kappa_{j\bi}}$
satisfying~$\kappa_{1\bar1}+\kappa_{2\bar2}+\kappa_{3\bar3}=0$
and
\be
\label{eq: G2streqs}
\begin{aligned}
\d\alpha_i &= -\kappa_{i\bar{l}}\w\alpha_l 
                + \overline{\strut\alpha_j\w\alpha_k}\\
\d\kappa_{i\bj} &= -\kappa_{i\bk}\w\kappa_{k\bj} 
                  + {\ts\frac34}\alpha_i\w\overline{\alpha_j}
                  -{\ts\frac14}\delta_{i\bj}\,\alpha_l\w\overline{\alpha_l}
\end{aligned}
\ee
where, in the first equation,~$(i,j,k)$ is an even permutation of~$(1,2,3)$
and where
\be
u^*\omega = {\ts\frac\iC2}\,\bigl(\alpha_1\w\overline{\alpha_1}+
    \alpha_2\w\overline{\alpha_2}+\alpha_3\w\overline{\alpha_3}\bigr)
\qquad\text{and}\qquad
u^*\psi = \alpha_1\w\alpha_2\w\alpha_3\,.
\ee
For details, see%
\footnote{Note, however, that, in order to make that notation match with
the notation in this article, one must take~$\alpha_i = -2\theta_i$.}
\cite[Proposition~2.3]{MR84h:53091}.  
In particular, it follows that
the almost complex structure~$J$ is strictly quasi-integrable and that
$\omega = \omega(J)$ and~$\psi = \psi(J)$.  

Note that the $\eusu(3)$-valued matrix~$\kappa = (\kappa_{i\bj})$
is a connection matrix for the tangent bundle~$E = T^{1,0}S^6$,
and that this is compatible with the given Hermitian structure on~$E$.
The second formula in~\eqref{eq: G2streqs} then shows that this 
connection defines a pseudo-holomorphic structure on~$E$ and, 
in agreement with Proposition~\ref{prop: qiimpliespdoHYM}, 
it is indeed pseudo-Hermitian-Yang-Mills.

Note also that because~$\pbar\omega=0$, it follows (see Remark~\ref{rem: Hodge})
that the closed~$(1,1)$-forms on~$S^6$ are also coclosed.
In particular, the group $Z^{1,1}_0(S^6)$ 
injects into the trivial deRham group~$H^2(S^6,\bbC)$ and hence is trivial.
Thus, the only global pseudo-holomorphic line bundles over~$S^6$ 
are trivial and, indeed, any pseudo-holomorphic structure on a unitary
bundle~$E$ of rank~$r$ must have the trace of its curvature vanishing, 
so that the holonomy of such a connection necessarily lies in~$\SU(r)$.

\subsection{The flag manifold~$\SU(3)/\bbT^2$}
\label{ssec: nKflagmfd}
As an example of a different sort, consider the flag manifold~$M
=\SU(3)/\bbT^2$, where~$\bbT^2\subset\SU(3)$ is the maximal torus
consisting of diagonal matrices.  The canonical left invariant form 
on~$\SU(3)$ can be written in the form
\be
\gamma = g^{-1}\,\d g
= \bpm \iC\,\beta_1 & \phm\alpha_3 & -\overline{\alpha_2}\strut\\
        -\overline{\alpha_3}\strut & \iC\,\beta_2 & \phm\alpha_1\\
        \phm\alpha_2 & -\overline{\alpha_1}\strut & \iC\,\beta_3 \epm
\ee
where the~$\beta_i$ are real and satisfy~$\beta_1 + \beta_2 + \beta_3 = 0$.
The structure equation~$\d\gamma = -\gamma\w\gamma$ is then equivalent
to
\be
\label{eq: SU3streqs}
\begin{aligned}
\d\alpha_i &= -\iC(\beta_j-\beta_k)\w\alpha_i
                 +\overline{\strut\alpha_j\w\alpha_k}\\
\d(\iC\beta_i) &= \alpha_k\w\overline{\alpha_k}
                   -\alpha_j\w\overline{\alpha_j}
\end{aligned}
\ee
where~$(i,j,k)$ is an even permutation of~$(1,2,3)$.  It follows that,
if~$\pi:\SU(3)\to\SU(3)/\bbT^2$ is the coset projection, then there
is a unique almost complex structure~$J$ on~$\SU(3)/\bbT^2$
such that
\be
\begin{aligned}
\pi^*\bigl(\omega(J)\bigr) 
&= {\ts\frac\iC2}\,\bigl(\alpha_1\w\overline{\alpha_1}
    +\alpha_2\w\overline{\alpha_2}+\alpha_3\w\overline{\alpha_3}\bigr)\\
\pi^*\bigl(\psi(J)\bigr) &= \alpha_1\w\alpha_2\w\alpha_3\,.
\end{aligned}
\ee
As the structure equations indicate, 
$\bigl(J,\omega(J)\bigr)$ is strictly quasi-integrable and satisfies
\be
\d\bigl(\omega(J)\bigr) = 3\Im\bigl(\psi(J)\bigr)
 \qquad\text{and}\qquad 
\d\bigl(\psi(J)\bigr) = 2\,\bigl(\omega(J)\bigr)^2.
\ee

In particular,~$\pbar\bigl(\omega(J)\bigr)=0$, so that, 
just as in the case of~$S^6$ (again via Remark~\ref{rem: Hodge}),
closed~$(1,1)$-forms are both primitive and coclosed.  Thus, 
$Z^{1,1}_0(M)$ injects into~$H^2(M,\bbC)\simeq\bbC^2$ via the
map that sends closed forms to their cohomology classes.  
By~\eqref{eq: SU3streqs}, there are closed~$(1,1)$-forms~$B_i$ 
on~$M$ that satisfy~$\pi^*(B_i) = \d\beta_i$.  These $(1,1)$-forms
satisfy~$B_1+B_2+B_3=0$, but are otherwise indepdendent as forms
and (hence) as cohomology classes.  In particular, it follows 
that every element of~$H^2(M,\bbZ)$ can be represented by a
closed~$(1,1)$-form.  Thus, every complex line bundle over~$M$
carries a pseudo-holomorphic structure and this is unique up to
gauge equivalence.  In particular, it follows that the complex
vector bundles generated by these line bundles (and their duals)
all carry pseudo-holomorphic structures.  

Note that these pseudo-holomorphic line bundles exist even though 
there are clearly no almost complex hypersurfaces (even locally) in~$M$.

\subsection{Nearly K\"ahler structures}
\label{ssec: nKstructs}
Both of the above examples fall 
into a larger class of elliptic quasi-integrable structures 
that has been studied extensively,
that of \emph{nearly K\"ahler structures}.  
For references on what is known about these structures, 
the reader might consult~\cite{MR0267502} or~\cite{MR1707644}.  

For the purposes of this paper, it will be
convenient to take a definition of `nearly K\"ahler' that
is well-adapted for use in dimension~$6$.  (This is not the 
standard definition, but see~\cite{MR1707644}.)  
An $\SU(3)$-structure on a $6$-manifold~$M$ can be specified
by giving its fundamental $(1,1)$-form~$\omega$ (which defines
the Hermitian structure) and a its fundamental $(3,0)$-form~$\psi$.
These forms are required to satisfy the equations
\be
\d\omega = 3c\Im(\psi)  \qquad\text{and}\qquad 
\d\psi = 2c\,\omega^2
\ee
for some real constant~$c$.  

When~$c=0$, one sees immediately that
the underlying almost complex structure is integrable, that $\omega$
defines a K\"ahler metric for this complex structure, and that~$\psi$
is a parallel holomorphic volume form on this K\"ahler manifold.
In other words, such an $\SU(3)$-structure is what is usually called
`Calabi-Yau'.  When~$c\not=0$ (in which case, one says that the 
structure is \emph{strictly nearly K\"ahler}), one can consider instead 
the scaled pair~$\bigl(c^2\omega, c^3\psi\bigr)$ (which will have the
same underlying almost complex structure, but a metric scaled by~$c^2$)
and see that one is reduced to the case~$c=1$.  

A calculation using the structure equations~\cite{MR1707644} shows that,
if~$\pi:\Fs\to M$ is the $\SU(3)$-bundle over~$M$ whose local sections are
the special unitary coframings~$\alpha:TU\to\bbC^3$, then the tautological
forms and connections on~$\Fs$ satisfy
\be
\pi^*\omega = {\ts\frac\iC2}\,\bigl(\alpha_1\w\overline{\alpha_1}+
    \alpha_2\w\overline{\alpha_2}+\alpha_3\w\overline{\alpha_3}\bigr)
\qquad\text{and}\qquad
\pi^*\psi = \alpha_1\w\alpha_2\w\alpha_3\,
\ee
and
\be
\label{eq: nKstreqs}
\begin{aligned}
\d\alpha_i &= -\kappa_{i\bar{l}}\w\alpha_l 
                + c\,\,\overline{\strut\alpha_j\w\alpha_k}\\
\d\kappa_{i\bj} &= -\kappa_{i\bk}\w\kappa_{k\bj} 
                  + c^2\,\bigl({\ts\frac34}\alpha_i\w\overline{\alpha_j}
                -{\ts\frac14}\delta_{i\bj}\,\alpha_l\w\overline{\alpha_l}\bigr)
                  + K_{i\bj p\bar q}\,\alpha_q\w\overline{\alpha_{p}}\,,
\end{aligned}
\ee
where, in the first line of~\eqref{eq: nKstreqs}, 
$(i,j,k)$ is any even permuation of~$(1,2,3)$
and where~$K_{i\bj p\bar q}=K_{p\bj i\bar q}=K_{i\bar q p\bj} 
= \overline{K_{j\bi q\bar p}}$ and~$K_{i\bi p\bar q} = 0$.

In particular, note that the $\eusu(3)$-valued connection form~$\kappa$
has its curvature of type~$(1,1)$.  Thus,~$\kappa$ defines
a pseudo-holomorphic structure on the tangent bundle of~$M$, one
that is, in fact, pseudo-Hermitian-Yang-Mills. 

Note also that the curvature tensor~$K$ that shows up in these structure
equations takes values in an $\SU(3)$-irreducible bundle of (real) rank~$27$
and that it vanishes exactly when the structure is either flat (if~$c=0$)
or, up to a constant scalar factor, equivalent to the $\Gtwo$-invariant
structure on~$S^6$ (if~$c\not=0$).  

In particular, the well-known
result that the underlying metric of a strictly nearly K\"ahler structure
on a $6$-manifold is Einstein (with positive Einstein constant) 
follows immediately from this.  It also follows that all such structures 
are real analytic in, say, coordinates harmonic for the metric.

A straightforward calculation using exterior differential systems
shows that, modulo diffeomorphism, the local nearly K\"ahler structures
in dimension~$6$ depend on two arbitrary functions of~$5$ variables,
i.e., they have the same local generality as the well-known~$c=0$ case.

Meanwhile, no compact example with~$c\not=0$
that is not homogeneous appears to be known.

\begin{proposition}[Supercriticality]
\label{prop: supercritical}
Any nearly K\"ahler structure is critical 
for all of the functionals~$F_c$ of Theorem~\ref{thm: InvFuncts}.
\end{proposition}

\begin{proof}
This follows immediately from the definitions 
and Proposition~\ref{prop: nrlypsdoKahler}.
\end{proof}

\subsection{Projectivized tangent bundles}
\label{ssec: projtanbundles}
Let~$S$ be a complex surface endowed with a K\"ahler form~$\eta$
and let~$\pi:F\to S$ be the associated~$\U(2)$-coframe bundle,
with structure equations
\be
\d\bpm \eta_1\\\eta_2\epm
 = - \bpm \psi_{1\bar1} & \psi_{1\bar2}\\ \psi_{2\bar1} & \psi_{2\bar2}\epm
     \w \bpm \eta_1\\\eta_2\epm
\ee
where~$\pi^*\eta = \frac\iC2\bigl(\eta_1\w\overline{\eta_1}+
\eta_1\w\overline{\eta_1}\bigr)$ and~$\overline{\psi_{i\bj}}=-\psi_{j\bi}$
and
\be
\d\psi_{i\bj} = -\psi_{i\bk}\w\psi_{k\bj} 
                + {\ts\frac12}K_{i\bj\ell\bk}\,\eta_k\w\overline{\eta_\ell}\,
\ee
where, as usual~$K_{i\bj\ell\bk} = K_{\ell\bj i\bk} 
=  K_{i\bk\ell\bj} = \overline{K_{j\bi k\bar\ell} }\,$. 

Setting~$\alpha_1 = \eta_1$, $\alpha_2 = \overline{\eta_2}$, 
and $\alpha_3 = -\psi_{2\bar1}$, one computes
\be
\left.
\begin{aligned}
\d\alpha_1 
 &\equiv \phantom{{\ts\frac12}K_{1\bar12\bar2}}
     \,\overline{\alpha_2}\w\overline{\alpha_3}\\
\d\alpha_2 
 &\equiv \phantom{{\ts\frac12}K_{1\bar12\bar2}}
      \,\overline{\alpha_3}\w\overline{\alpha_1}\\
\d\alpha_3 
 &\equiv {\ts\frac12}K_{1\bar12\bar2}
       \,\overline{\alpha_1}\w\overline{\alpha_2}\\
\end{aligned} \quad
\right\}
\mod \alpha_1,\,\alpha_2,\,\alpha_3\,.
\ee
Letting~$\bbT^2\subset\U(2)$ be the subgroup of diagonal unitary matrices,
it follows that there is a well-defined almost complex structure~$J$
on~$M = F/\bbT^2$ whose $(1,0)$-forms pull back to~$F$ to be linear
combinations of~$\alpha_1$, $\alpha_2$, and~$\alpha_3$.  

Since, after pullback to~$F$, one has the formula,
\be
\omega(J) = {\ts\frac\iC2}\,
\left({\ts\frac12}K_{1\bar12\bar2}\,\alpha_1\w\overline{\alpha_1}
    + {\ts\frac12}K_{1\bar12\bar2}\,\alpha_2\w\overline{\alpha_2}
    + \alpha_3\w\overline{\alpha_3}\right),
\ee
it follows that~$\omega(J)$ will be nondegenerate on~$M$ (which
is the projectivized tangent bundle of~$S$) as long as the
K\"ahler metric~$\eta$ on~$S$ 
has nonvanishing holomorphic bisectional curvature.
In particular, $J$ is elliptically quasi-integrable if~$\eta$
has positive holomorphic bisectional curvature and hyperbolically
quasi-integrable if~$\eta$ has negative holomorphic bisectional
curvature.  In any case, the Nijenhuis tensor of~$J$ has real type.

Using the structure equations on~$F$, one also sees that, if~$\eta$
has constant holomorphic bisectional curvature~$c$, then the almost
complex structure~$J$ will be supercritical.  When~$c>0$, this is
just the nearly K\"ahler case of~$\SU(3)/\bbT^2$ already discussed
and is quasi-integrable. When~$c<0$, i.e., when~$(S,\eta)$ is a 
ball quotient (i.e., locally isometric to the Hermitian symmetric 
space~$\SU(1,2)/\U(2) = \bbB^2$), 
this gives examples (some of which are compact)
of supercritical almost complex structures 
that are hyperbolically quasi-integrable. 

\begin{remark}[Hyperbolic nearly K\"ahler structures]
This last example points out that there is a hyperbolic
analog of nearly K\"ahler structures that is worth defining:
An $\SU(1,2)$-structure on~$M^6$, with defining $(1,1)$-form~$\omega$
of Hermitian signature~$(1,2)$ and $(3,0)$-form~$\psi$ 
will be said to be \emph{hyperbolically nearly K\"ahler} if
it satisfies the structure equations
\be
\d\omega = 3c\,\Im(\psi)\qquad\text{and}\qquad
\d\psi = 2c\,\omega^2
\ee
for some (real) constant~$c$.  Again, the usual methods of exterior
differential systems shows that the local (i.e., germs of) 
hyperbolic nearly K\"ahler structures modulo diffeomorphism 
depend on $2$ arbitrary functions of~$5$ variables, just
as in the elliptic case.  
\end{remark}

\subsection{Generality}
\label{ssec: generality}
It is natural to ask how `general' (local) quasi-integrable  
almost complex structures are.  This is something of a vague
question, but the notion of generality can be made precise in a 
number of ways, usually (when considering the generality of a 
space of solutions of a system of PDE) by an appeal to Cartan's 
theory of systems in involution~\cite{MR1083148}. 

For the sake of simplicity, 
I will only consider the elliptic quasi-integrable case
and only occasionally remark on the analogous results in
the hyperbolic case.

Since quasi-integrability is a first-order condition on~$J$,
the set of $1$-jets of quasi-integrable structures is a 
well-defined subset~$\QJ(M)\subset J^1\bigl(\Js(M)\bigr)$.  
Unfortunately, $\QJ(M)$ is not smooth, but, at its smooth points,
it has codimension~$8$ in~$J^1\bigl(\Js(M)\bigr)$. 

By contrast, quasi-integrability for~$\U(3)$-structures 
is better behaved as a first-order \textsc{pde} system:  
Let~$\Us(M)\to M$ be the bundle of unitary structures over~$M$,
whose general fiber is modeled on~$\GL(6,\bbR)/\U(3)$.
Then the set of~$1$-jets of quasi-integrable~$\U(3)$-structures
on~$M$ is a smooth submanifold~$\QU(M)\subset J^1\bigl(\Us(M)\bigr)$
of codimension~$16$.  

Unfortunately, the condition of quasi-integrability 
for~$\U(3)$-structures turns out not to be involutive, 
and a full analysis has yet to be done.  The point of this 
subsection is to explain at least where the first 
obstruction appears.  Along the way, some useful information 
about quasi-integrable $\U(3)$-structures will be uncovered.

\subsubsection{The structure equations}
Since quasi-integrability for a $\U(3)$-structure~$(J,\eta)$ 
is well-defined as a condition 
on the underlying~$\bbR^+{\cdot}\U(3)$-structure, 
it is natural to study the structure equations 
on the associated~$\bbR^+{\cdot}\U(3)$-bundle~$\pi:\Fs\to M$
whose local sections are the `conformally unitary' coframings
$\alpha:TU\to\bbC^3$, i.e., $\alpha$ is~$J$-linear and satisfies
${{\eta\vrule}_{U}} = \mu\,\frac\iC2\,(\alpha^1\w\overline{\alpha^1}
+ \alpha^2\w\overline{\alpha^2} + \alpha^3\w\overline{\alpha^3})$
for some positive function~$\mu$.

A straightforward analysis of the intrinsic torsion
of the~$\bbR^+{\cdot}\U(3)$-structure~$\Fs$ 
shows that~$\zeta = (\zeta_i)$, 
the canonical~$\bbC^3$-valued $1$-form on~$\Fs$, 
satisfies unique first structure equations of the form
\be
\label{eq: dzeta}
\d\zeta_i 
= -\rho\w\zeta_i - \kappa_{i\bj}\w\zeta_j
  - B_{ij\bk}\,\overline{\zeta_j}\w\zeta_k
  + {\ts\frac12}\varepsilon_{ijk}\,\lambda\,\overline{\zeta_j\w\zeta_k}
\ee
where~$\rho=\overline{\rho}$ is real, 
$\kappa_{i\bj} = -\overline{\kappa_{j\bi}}$,
and~$B_{i\ell\bj}=-B_{\ell i\bj}$ satisfying~$B_{ij\bj}=0$ 
and~$\lambda$ are complex functions defined on~$\Fs$.
(N.B.: The unitary summation convention 
is being employed and, as usual, the symbol~$\varepsilon_{ijk}$ 
is completely antisymmetric and satisfies~$\varepsilon_{123}=1$.)
The $1$-forms~$\rho$ and~$\kappa = (\kappa_{i\bj})=-{}^t\bar\kappa$ 
are the~$\bbR$- and~$\euu(3)$-components of the canonical connection form 
on the $\bbR^+{\cdot}\U(3)$-bundle~$\pi:\Fs\to M$.  The
functions~$\lambda$ and~$B_{ij\bk}$ represent the components of
the first-order tensorial (torsion) invariants 
of the $\bbR^+{\cdot}\U(3)$-structure.

\begin{remark}[Relation with Gray-Hervella torsion]
\label{rem: GrayHervellatorsion}
The reader who is familiar with the canonical connection 
of a~$\U(3)$-structure and its four irreducible torsion components
as described by Gray and Hervella in~\cite{MR0581924} 
may be wondering how these structure equations correspond.  
The function~$\lambda$ represents their tensor~$W_1$, 
the quasi-integrable assumption is the vanishing of their tensor~$W_2$, 
and the functions~$B_{ij\bk}$ represent the components of their tensor~$W_3$.  
Their tensor~$W_4$ is not conformally invariant and, 
as a result, does not appear in the structure equations~\eqref{eq: dzeta}.
Instead, one has the uniqueness of the conformal connection form~$\rho$,
corresponding to the fact that the first prolongation of the Lie algebra
of~$\bbR^+{\cdot}\U(3)\subset\GL(6,\bbR)$ is trivial.
\end{remark}

Note that, for the underlying almost complex structure~$J$, 
the canonical forms~$\omega$ and~$\psi$ satisfy
\be
\pi^*\omega = {\ts\frac\iC2}|\lambda|^2\,\zeta_i\w\overline{\zeta_i}
\qquad\text{and}\quad
\pi^*\psi = |\lambda|^2\overline{\lambda}\,\,\zeta_1\w\zeta_2\w\zeta_3\,.
\ee

The second structure equations are computed by taking the exterior
derivative of the first structure equations.  In order to do this
most efficiently, first set
\be
\begin{aligned}
D\lambda &= \d\lambda - \lambda \bigl(\rho - \tr\kappa\bigr),\\
DB_{ij\bk}  &= \d B_{ij\bk} 
             + B_{pj\bk}\kappa_{i\bar p} + B_{ip\bk}\kappa_{j\bar p}
             - B_{ij\bar p}\kappa_{p\bk} - B_{ij\bk}\,\rho,\\
K_{i\bj} &=  \d\kappa_{i\bj} + \kappa_{i\bk}\w\kappa_{k\bj}
= -\overline{K_{j\bi}}.
\end{aligned}
\ee

Then the exterior derivative of~\eqref{eq: dzeta} becomes, after
some rearrangment and term collection,%
\footnote{Although~$\overline{\varepsilon_{ijk}}=\varepsilon_{ijk}$, 
I use the former when needed to maintain the unitary summation convention
that summation is implied when an index appears both barred and unbarred
in a given term.  N.B.:  In an expression like~$\overline{a_\bi}$, the
index~$i$ is to be regarded as unbarred.}
\be
\label{eq: ddzeta}
\begin{aligned}
0 & = {}-\d\rho\w\zeta_i - K_{i\bj}\w\zeta_j          
         - DB_{ij\bk}\w\overline{\zeta_j}\w\zeta_k
      + {\ts\frac12}\varepsilon_{ijk}\,
          D\lambda\w\overline{\zeta_j}\w\overline{\zeta_k}\\
& \quad\quad\quad
  + \bigl(B_{i\ell\bar p}\overline{B_{\ell q \bj}}
          - {\ts\frac12}\varepsilon_{ijk}\overline{\varepsilon_{pqk}}\,
              |\lambda|^2\bigr)\, \zeta_p\w\zeta_q\w\overline{\zeta_j}\\
&\qquad\qquad 
   +\bigl(B_{iq\bj}B_{j p \bk}
          +\varepsilon_{iq\ell}\lambda\,\overline{B_{\ell k\bar p}}\,\bigr)
            \,\zeta_k\w\overline{\zeta_p}\w\overline{\zeta_q}\,.
\end{aligned}
\ee

By exterior algebra, there exist unique functions~$\lambda_\ell$,
$\lambda_{\bar\ell}$, $B_{ij\bk\bar\ell}$,~$B_{ij\bk\ell}$,~$S_{ij}=-S_{ji}$,
$R_{i\bj}= \overline{R_{j\bi}}$, $L_{i\bj k\ell}=-L_{i\bj\ell k}$, 
and $M_{i\bj k\bar\ell}=\overline{M_{j\bi\ell\bk}}$ on~$\Fs$ such that
\be
\label{eq: dDlambdaDBdrhodkappa}
\begin{aligned}
D\lambda &=\lambda_{\bar\ell}\,\zeta_\ell+\lambda_\ell\,\overline{\zeta_\ell}\,,\\
DB_{ij\bk}
&= B_{ij\bk\bar\ell}\,\zeta_\ell+B_{ij\bk\ell}\,\overline{\zeta_\ell}\,,\\
\d\rho & = {\ts\frac12}S_{k\ell}\,\overline{\zeta_k}\w\overline{\zeta_\ell}
           + {\ts\frac\iC2}\,R_{k\bj}\,\zeta_j\w\overline{\zeta_k}
           + {\ts\frac12}\overline{S_{k\ell}}\,\zeta_k\w\zeta_\ell\,,\\
K_{i\bj} 
&= {\ts\frac12}L_{i\bj k\ell}\,\overline{\zeta_k}\w\overline{\zeta_\ell}
      + {\ts\frac12}\,M_{i\bj k\bar\ell}\,\zeta_\ell\w\overline{\zeta_k}
      - {\ts\frac12}\overline{L_{j\bi k\ell}}\,\zeta_k\w\zeta_\ell\,.
\end{aligned}
\ee
The coefficients appearing in the terms on the right hand side 
of~\eqref{eq: dDlambdaDBdrhodkappa} are the second-order invariants
of the $\bbR^+{\cdot}\U(3)$-structure~$\Fs$.  

Substituting these equations into~\eqref{eq: ddzeta} 
yields the \emph{first Bianchi identities}, which can be sorted
into~$(p,q)$-types in the obvious way.  I now want to consider some of
the consequences of these identities.

First, let~$\sigma:\bbR^+{\cdot}\U(3)\to\bbC^*$ be the representation that
satisfies~$\sigma(tg) = t\det(\bar g)$ for~$t>0$ and~$g\in\U(3)$
and define the complex line bundle~$\Lambda = \Fs\times_\sigma\bbC$ over~$M$.
Then the $1$-form~$\rho-\tr\kappa$ can be viewed 
as defining a connection~$\nabla$ on the line bundle~$\Lambda$
and the equation
\be
\d\lambda = \lambda(\rho-\tr\kappa) 
          + \lambda_{\bar\ell}\,\zeta_\ell+\lambda_\ell\,\overline{\zeta_\ell}
\ee
shows that~$\lambda$ represents a smooth section~$L$ of~$\Lambda$.
Consequently, the $1$-forms~$\lambda_\bi\,\zeta_i$ 
and~$\lambda_i\,\overline{\zeta_i}$ represent, 
respectively~$\nabla^{1,0}L$ and~$\nabla^{0,1}L = \pbar L$.

\begin{proposition}[Holomorphicity of~$L$]
\label{prop: Lholomorphic}
Let~$M^6$ be endowed with an~$\bbR^+{\cdot}\U(3)$-structure
that is quasi-integrable.
Then the section~$L$ is pseudo-holomorphic, i.e., $\pbar L = 0$.
In particular, if~$M$ is connected, then either~$L$ vanishes identically
\oppar in which case the almost complex structure~$J$ is 
integrable\clpar\ or else the zero locus of~$L$ 
\oppar and hence, of the Nijenhuis 
tensor\clpar\ is of \oppar real\clpar\ codimension at least~$2$ in~$M$.
\end{proposition}

\begin{proof}
One immediately sees that the~$(0,3)$-component 
of the first Bianchi identities is the three equations
\be
\lambda_\ell = 0.
\ee
Thus~$\pbar L=0$, as claimed.

The remainder of the proposition now follows by a standard argument 
about pseudo-holomorphic sections of line bundles
endowed with connections over almost complex manifolds, 
so I will only sketch it.  

First, note that, for any connected pseudo-holomorphic 
curve~$C\subset M$, the restriction of~$(\Lambda,\nabla)$ to~$C$ 
becomes a complex line bundle~$\Lambda_C$ that carries a unique holomorphic
structure for which~$\nabla^{0,1} = \pbar_C$ and that the restriction of
the section~$L$ to~$C$ becomes a holomorphic section of~$\Lambda_C$ with
respect to this holomorphic structure.  In particular, $L$ either vanishes
identically on~$C$ or else has only isolated zeros of finite order.  

Consequently, if~$L$ vanishes to infinite order at any point~$p\in M$, 
it must vanish identically on any connected pseudo-holomorphic curve~$C$ in~$M$
that passes through~$p$.  Since the union of the pseudo-holomorphic discs in~$M$
passing through~$p$ contains an open neighborhood of~$p$,
it follows that the set of points~$Z_\infty(L)$ at which~$L$ 
vanishes to infinite order is an open set.  Since~$Z_\infty(L)$ 
is clearly closed and~$M$ is connected, it follows that either~$Z_\infty(L)$
is all of~$M$, so that~$L$ vanishes identically, or else that~$Z_\infty(L)$ 
is empty.

Now suppose that~$Z_\infty(L)$ is empty.   If~$Z(L)$, the set of points
in~$M$ at which~$L$ vanishes is empty, then there is nothing to prove,
so suppose that~$L(p)=0$.  Since~$L$ vanishes to finite order at~$p$,
there is an immersed pseudo-holomorphic 
disk~$D\subset M$ passing through~$p$ such that~$L$ does not vanish identically 
on~$D$.  Since~$L$ restricted to~$D$ is a holomorphic section of a holomorphic 
line bundle, it follows that, by shrinking~$D$ if necessary, it can be assumed
that~$p$ is the only zero of~$L$ on~$D$ and that~$\nu_p(L_D)=k$ for
some integer~$k>0$.

One can now embed~$D$ into a (real) $4$-parameter
family of pseudo-holomorphic disks~$D_t\subset M$
with~$t\in\bbR^4$ that foliate a neighborhood of~$p$
and satisfy~$D_0 = D$.  One can further assume that~$L$ does not vanish
on any of the boundaries~$\p D_t$ (since it does not vanish on~$\p D_0$).
Then for each~$t$, the section~$L$ restricted to~$D_t$ has a zero locus 
consisting of~$k$ points (counted with multiplicity).
Thus, the intersection~$Z(L)\cap D_t$ is a finite set of points 
whose total multiplicity is~$k$.  That the real codimension of~$Z(L)$ 
is at least~$2$ now follows from this.
\end{proof}

\begin{remark}[Speculation about~$Z(L)$]
It is tempting to conjecture that if~$Z(L)$ 
is neither empty nor all of~$M$, then it has a stratification
\be
Z(L) = Z^1(L)\cup Z^2(L)\cup Z^3(L)
\ee
where~$Z^i(L)\subset M$ is a smooth almost complex
submanifold of~$(M,J)$ of complex codimension~$i$, 
with~$Z(L) = \overline{Z^1(L)}$.  One might even expect
that, when~$M$ is compact,~$Z^1(L)$, counted with the appropriate
multiplicity, would define a class in~$H_4(M,\bbZ)$ that would
be dual to~$c_1(\Lambda)$.

In fact, let~$\nu:Z(L)\to\bbZ^+$ be the upper semicontinuous function 
giving the order of vanishing of~$L$.  Then the set~$Z^*(L)\subset Z(L)$
on which~$\nu$ is locally constant is open and dense in~$Z(L)$ 
and it is not difficult to show that~$Z^*(L)$ is a smooth $4$-dimensional 
submanifold of~$M$ whose tangent spaces are everywhere complex.
The difficulty apprears to be in understanding the `singular' part
of~$L$, i.e., the set of points at which~$\nu$ is not locally constant.

Note that, in the open dense set~$M^* = M\setminus Z(L)$, 
the almost complex structure~$J$ is elliptically strictly 
quasi-integrable and hence there cannot be any complex 
codimension~$1$ almost complex submanifolds in~$M^*$, even locally.
\end{remark}

\begin{remark}[Singularity issues]
\label{rem: sing}
Let~$S$ be a complex surface with a K\"ahler form~$\eta$
whose bisectional curvature is everywhere nonnegative.  Then,
as explained in~\S\ref{ssec: projtanbundles}, the projectivized
tangent bundle~$M = F/\bbT^2$ carries an almost complex 
structure~$J$ whose Nijenhuis tensor is everywhere of real type.
On the open set~$M^*\subset M$ that represents the points in the
projectivized tangent bundle where the holomorphic bisectional curvature
is strictly positive, $J$ is strictly elliptically quasi-integrable.
However, the set of points in~$M$ where the holomorphic
bisectional curvature (represented by~$K_{1\bar12\bar2}\ge0$)
vanishes can be quite wild.  In particular, it does not have
to have real codimension~$2$ in~$M$.  For example, it can
easily happen that this zero locus is a single point in~$M$.
By Proposition~\ref{prop: Lholomorphic},
such elliptically quasi-integrable almost complex structures
do not underlie an elliptically quasi-integrable unitary structure on~$M$.  
This example gives an indication of what sort of singularities
are being avoided by working with quasi-integrable unitary structures
rather quasi-integrable almost complex structures.
\end{remark}

\begin{remark}[Second-order forms]
\label{eq: 2ndorderforms}
Note that a consequence of the holomorphicity of~$L$
is the following formulae:
First,~$\d\psi = \pbar\psi + 2\,\omega^2$, where
\be
\pi^*(\pbar\psi) 
 = 2|\lambda|^2\,\overline{\bigl(\lambda_{\bk}\,\zeta_k\bigr)}
                \w (\zeta_1\w\zeta_2\w\zeta_3).
\ee
Second~$\d\omega = 3\Im(\psi) +\p\omega + \pbar\omega$, where
\be
\pi^*(\pbar\omega) 
  = {\ts\frac\iC2}\lambda\,\overline{\bigl(\lambda_{\bk}\,\zeta_k\bigr)}
                \w (\zeta_\ell\w\overline{\zeta_\ell})
    -{\ts\frac\iC2}|\lambda|^2 B_{ij\bk}\,\zeta_k\w\overline{\zeta_i\w\zeta_j}.
\ee
The two terms on the right hand side 
represent the decomposition of~$\pbar\omega$ 
into its primitive types, as defined 
by the underlying~$\bbR^+{\cdot}\U(3)$-structure.
\end{remark}

I will now resume the investigation of the first Bianchi identities.
The $(3,0)$-part of the first Bianchi identity becomes
\be
\label{eq: 1stB30}
\overline{\varepsilon_{ijk}}\,S_{jk}
- \overline{\varepsilon_{jk\ell}}\,L_{j\bi k\ell} = 0,
\ee
which shows that~$S_{jk}$ is determined in terms of~$L_{j\bi k\ell}$.
The $(2,1)$-part of the first Bianchi identity is the more complicated
\be
\label{eq: 1stB21}
\begin{aligned}
{\ts\frac\iC2}\bigl(\delta_{i\bj}R_{k\bar\ell}-\delta_{i\bar\ell}R_{k\bj}\bigr)
+{\ts\frac12}\bigl(M_{i\bj k\bar\ell}-M_{i\bar\ell k\bj}\bigr)
+ B_{ik\bj\bar\ell}-B_{ik\bar\ell\bj}\qquad & \\
{} + B_{ip\bar\ell}\overline{B_{pj\bk}}-B_{ip\bj}\overline{B_{p\ell\bk}}
+\varepsilon_{ikp}\overline{\varepsilon_{j\ell p}}\,|\lambda|^2 &=0.
\end{aligned}
\ee
The~$(1,2)$-part of the first Bianchi identity is
\be
\label{eq: 1stB12}
\begin{aligned}
\delta_{i\bar\ell} S_{jk}+L_{i\bar\ell jk}
+ \bigl(B_{ik\bar\ell j}{-}B_{ij\bar\ell k}\bigr)
-\varepsilon_{ijk}\lambda_{\bar\ell}\qquad\qquad\qquad & \\
{}+ B_{ij\bar p}B_{pk\bar\ell}-B_{ik\bar p}B_{pj\bar\ell}
+\bigl(\varepsilon_{ijp}\overline{B_{p\ell\bk}}
          -\varepsilon_{ikp}\overline{B_{p\ell\bj}}\,\bigr)\,\lambda &=0.
\end{aligned}
\ee
Finally, remember that the $(0,3)$-part of the Bianchi identity has already
been seen to be
\be
\label{eq: 1stB03}
\lambda_\ell = 0.
\ee

For what I have in mind, it will be necessary to solve these equations 
more-or-less explicitly in order to understand how many free derivatives
there are at second-order for quasi-integrable structures.  After some
exterior algebra and representation theoretic arguments, one finds that
the following method works:

It is useful to introduce notation for the traces of the covariant
derivatives of the~$B$-tensor.  Thus, set
\be
B_{ip\bj\bar p}  = u_{i\bj} + \iC\, v_{i\bj}
\ee
where~$u = (u_{i\bj})={}^t\bar u$ and~$v = (v_{i\bj})={}^t\bar v$
take values in Hermitian $3$-by-$3$ matrices.  Note that, 
because~$B_{ji\bj}=0$ by definition, it follows that~$u$ and~$v$ 
are themselves traceless.  Next, write
\be
B_{ij\bk k} = -B_{ji\bk k} =  4\varepsilon_{ijp}\,a_{\bar p}
\ee
for functions~$a = (a_{\bar p})$.  Then one finds that
\be
B_{ij\bk\ell} = B^*_{ij\bk\ell} 
                + \varepsilon_{ijp}\bigl(
                     \delta_{\ell\bar p}\,a_{\bk}
                   + \delta_{\ell\bk}\,a_{\bar p}\bigr)
\ee
where~$B^*_{ij\bk\ell}=-B^*_{ji\bk\ell}$ satisfies~$B^*_{ij\bk k}
= B^*_{ik\bk\ell} = 0$.  

With these quantities in hand, one finds that the Bianchi identities
are equivalent to the following relations:
\be
\label{eq: S1stBianchi}
S_{ij} = \varepsilon_{ijk}\bigl(4\,a_\bk + {\ts\frac32}\lambda_\bk\bigr),
\ee
\be
\label{eq: L1stBianchi}
L_{i\bj k\ell}
= B^*_{ik\bj\ell} - B^*_{i\ell\bj k}
   +\delta_{i\bj}\,Q_{k\ell}+\delta_{k\bj}\,P_{i\ell}-\delta_{\ell\bj}\,P_{ik}
   +\delta_{k\bj}\,F_{i\ell}-\delta_{\ell\bj}\,F_{ik}\,,
\ee
where
\be
\label{eq: PQF1stBianchi}
\begin{aligned}
Q_{ij} &= -\varepsilon_{ijk}\bigl(2\,a_\bk + {\ts\frac12}\lambda_\bk\bigr)
         = -Q_{ji}\\
P_{ij} &= -\varepsilon_{ijk}\bigl(3\,a_\bk 
                   +\phantom{\ts\frac12}\lambda_\bk\bigr)
         = - P_{ji}\\
F_{ij} &= {\ts\frac12}B_{ip\bar q}B_{jq\bar p} 
             + {\ts\frac12}\varepsilon_{pqi}\,\lambda\,\overline{B_{pq\bj}}
        = F_{ji}\,,
\end{aligned}
\ee
\be
\label{eq: F1stBianchi}
R_{i\bj} = -2\,v_{i\bj}
\ee
(in particular, note that~$R$ is traceless),
and, finally,
\be
\label{eq: M1stBianchi}
\begin{aligned}
M_{i\bj k\bar\ell} &= K_{i\bj k\bar\ell} 
      - \delta_{i\bj}\,u_{k\bar\ell} - \delta_{k\bar\ell}\,u_{i\bj}
   -\iC\bigl(\delta_{i\bar\ell}\,v_{k\bj}-\delta_{k\bj}\,v_{i\bar\ell}\bigr)\\
&\qquad + 2\delta_{i\bj}\,B_{pq\bar\ell}\overline{B_{pqk}}
        + 2\delta_{k\bar\ell}\,B_{pq\bj}\overline{B_{pq\bi}}
        -\delta_{i\bj}\delta_{k\bar\ell}
         \bigl(B_{pq\bar r}\overline{B_{pq\bar r}}+2|\lambda|^2\bigr)\,,
\end{aligned}
\ee
where~$K_{i\bj k\bar\ell} = K_{k\bj i\bar\ell} = K_{i\bar\ell k\bj}
= \overline{K_{j\bi\ell\bk}}$ has the symmetries of a K\"ahler curvature
tensor.

The main use of these formule will be in the following observation:
\emph{All of the second-order invariants of a~$\bbR^+{\cdot}\U(3)$-structure
are expressed in terms of the covariant derivatives of its $\lambda$-
and~$B$-tensors plus the `K\"ahler-component' of the curvature of its
$\kappa$-connection.  Moreover, the only relations among these 
latter three invariants, other than the usual symmetries of the K\"ahler curvature, are the relations~$\lambda_i=0$.} 

\begin{proposition}[Noninvolutivity of the quasi-integrable condition]
The first order system on~$\bbR^+{\cdot}\U(3)$-structures that 
defines quasi-integrability is not formally integrable.  
In particular, it is not involutive.
\end{proposition}

\begin{proof}
The structure equations and Bianchi identities imply the formulae
\be
\label{eq: drho1stBianchi}
\d\rho = 
\varepsilon_{k\ell p}\bigl(2\,a_{\bar p}+{\ts\frac34}\lambda_{\bar p}\bigr)\,
\overline{\zeta_k}\w\overline{\zeta_\ell} 
 - \iC\,v_{i\bj}\,\zeta_j\w\overline{\zeta_i}
+ \overline{\varepsilon_{k\ell p}\bigl(2\,a_{\bar p}
   +{\ts\frac34}\lambda_{\bar p}\bigr)}\,\zeta_k\w\zeta_\ell
\ee
and, by~\eqref{eq: L1stBianchi} and~\eqref{eq: PQF1stBianchi},
\be
\label{eq: dtrkappa1stBianchi}
\begin{aligned}
\d(\tr\kappa) &=
{\ts\frac12}L_{i\bi k\ell}\,\overline{\zeta_k}\w\overline{\zeta_\ell}
      + {\ts\frac12}\,M_{i\bi k\bar\ell}\,\zeta_\ell\w\overline{\zeta_k}
      - {\ts\frac12}\overline{L_{i\bi k\ell}}\,\zeta_k\w\zeta_\ell\\
&=-\varepsilon_{k\ell p}\bigl(6\,a_{\bar p}
      +{\ts\frac74}\lambda_{\bar p}\bigr)\,
\overline{\zeta_k}\w\overline{\zeta_\ell} 
 + {\ts\frac12}M_{i\bi k\bar\ell}\,\zeta_\ell\w\overline{\zeta_k}\\
&\qquad\qquad\qquad+ \overline{\varepsilon_{k\ell p}\bigl(6\,a_{\bar p}
                +{\ts\frac74}\lambda_{\bar p}\bigr)\,}\,\zeta_k\w\zeta_\ell
\end{aligned}
\ee

Meanwhile, the exterior derivative of the equation
\be
\label{eq: dlambdahol}
\d\lambda = \lambda\,\bigl(\rho - \tr\kappa\bigr) + \lambda_\bi\,\zeta_i 
\ee
yields
\be
0 = (\lambda_\bi\,\zeta_i )\w\bigl(\rho - \tr\kappa\bigr)
   + \lambda\,\bigl(\d\rho - \d(\tr\kappa)\bigr)
   + \d\lambda_{\bi}\w\zeta_i +\lambda_{\bi}\,\d\zeta_i\,
\ee
which, by the first structure equations, can be written in the form
\be
\label{eq: ddlambda}
0 =  \lambda\,\bigl(\d\rho - \d(\tr\kappa)\bigr)
     + D\lambda_{\bi}\w\zeta_i 
+{\ts\frac12}\varepsilon_{ijk}\,\lambda\lambda_{\bi}\,\overline{\zeta_j\w\zeta_k}\,,
\ee
where
\be
D\lambda_{\bi} 
= \d\lambda_{\bi} - (2\rho-\tr\kappa)\lambda_{\bi} -\kappa_{j\bi}\,\lambda_{\bj}
       - B_{j\ell\bi}\lambda_{\bj}\,\overline{\zeta_\ell}\,
\ee
is $\pi$-semi-basic.  Thus, taking the $(0,2)$-part of~\eqref{eq: ddlambda}
yields
\be
\label{eq: drho-dtrkappa}
\lambda\bigl(\d\rho - \d(\tr\kappa)\bigr)^{0,2} 
= -{\ts\frac12}\varepsilon_{ijk}\,\lambda\lambda_{\bi}\,
\overline{\zeta_j\w\zeta_k}\,.
\ee
On the other hand, \eqref{eq: drho1stBianchi} 
and \eqref{eq: dtrkappa1stBianchi} yield
\be
\label{eq: drho-dtrkappa1stBianchi}
\lambda\bigl(\d\rho - \d(\tr\kappa)\bigr)^{0,2} 
= \varepsilon_{ijk}\,\lambda
      \bigl(8\,a_\bi + {\ts\frac52}\lambda_\bi\bigr)\,
\overline{\zeta_j\w\zeta_k}\,.
\ee
In other words, the equation
\be
\label{eq: FIobstruction}
\lambda\bigl(8\,a_\bi + 3\lambda_\bi\bigr) = 0
\ee
is an identity for quasi-integrable~$\bbR^+{\cdot}\U(3)$-structures.

It is important to understand how this identity was derived.
While it is a relation on the second-order invariants, 
it is \emph{not} an algebraic consequence of the first Bianchi identity.  
Instead, it was derived by combining the first Bianchi identity 
with the result of differentiating the relation~\eqref{eq: dlambdahol},
which is, itself, a part of the first Bianchi identity.
It is the presence of this second-order identity 
that \emph{cannot} be found by differentiating only one time 
the first-order defining equations 
for quasi-integrable~$\bbR^+{\cdot}\U(3)$-structures 
that shows that these defining equations are not formally integrable 
and, hence, not involutive.  
(See Remark~\ref{rem: formalint} for a further discussion of this point.)
\end{proof}

\begin{corollary}[Quasi-integrable dichotomy]
\label{cor: qidichotomy}
Let~$M^6$ be connected and let~$\bigl(J,[\eta]\bigr)$
be a quasi-integrable $\bbR^+{\cdot}\U(3)$-structure on~$M$.
Then either~$J$ is integrable \oppar i.e., $\lambda$
vanishes identically\clpar\ or else the structure satisfies
the second-order equations
\be
\label{eq: QI2ndorderreln}
\lambda_{\bar\ell} 
= -{\ts\frac13}\overline{\varepsilon_{ij\ell}}\,B_{ij\bk k}\,.
\ee
\end{corollary}

\begin{proof}
By Proposition~\ref{prop: Lholomorphic}, 
either~$\lambda$ vanishes identically or else 
its zero locus has no interior.  In this latter case,
\eqref{eq: FIobstruction} implies that~$8\,a_\bi + 3\lambda_\bi$
must vanish identically.  Tracing through the definitions, 
this is~\eqref{eq: QI2ndorderreln}.
\end{proof}

\begin{remark}[Formal integrability]
\label{rem: formalint}
An adequate discussion of formal integrability would be too long
to include here; see~\cite[Ch.~IX]{MR1083148} for details.
Very roughly speaking, a system~$\Sigma$ of first-order \textsc{pde} 
is formally integrable if any differential equation 
of order~$q$ that is satisfied by all of the solutions 
of~$\Sigma$ is derivable from~$\Sigma$ by differentiating
the equations in~$\Sigma$ at most~$q{-}1$ times.

However, it is worth pointing out explicitly what is going
on in this particular case.  Let~$\pi:\CU(M)\to M$ denote
the bundle whose local sections are the local 
$\bbR^+{\cdot}\U(3)$-structures on~$M$. (The notation~$\CU(M)$
is meant to denote `conformal unitary'.)  This bundle has
fibers isomorphic to the $26$-dimensional homogeneous 
space~$\GL(6,\bbR)/\bigl(\bbR^+{\cdot}\U(3)\bigr)$.  (Thus,
such local structures depend on $26$ functions of $6$ variables.)

The $1$-jets of local quasi-integrable~$\bbR^+{\cdot}\U(3)$-structures
form a smooth subbundle~$Q\subset J^1\bigl(\CU(M)\bigr)$, 
one that has (real) codimension~$16$ in~$J^1\bigl(\CU(M)\bigr)$.
Its first prolongation~$Q^{(1)}\subset J^2\bigl(\CU(M)\bigr)$
is also a smooth submanifold and the projection~$Q^{(1)}\to Q$
is a smooth submersion.  However, 
its second prolongation~$Q^{(2)} = \bigl(Q^{(1)}\bigr)^{(1)}
\subset J^3\bigl(\CU(M)\bigr)$ is \emph{not} a smooth manifold
and its projection~$Q^{(2)}\to Q^{(1)}$ is not surjective
and not a submersion.

In fact, the relations~\eqref{eq: FIobstruction} show that
the image of the projection of~$Q^{(2)}\to Q^{(1)}$ has
real codimension at least%
\footnote{It has not been shown that~\eqref{eq: FIobstruction}
defines the image of~$Q^{(2)}$ in~$Q^{(1)}$, 
only that $Q^{(2)}$ lies in union of the two
submanifolds~$Q'_1\subset Q^{(1)}$ defined by~$\lambda=0$
(and hence of codimension~$2$) and~$Q''_1\subset Q^{(1)}$
defined by~$8\,a_\bi + 3\lambda_\bi=0$ (and hence of codimension~$6$).
While I believe that the union $Q'_1\cup Q''_1$ is the image
of~$Q^{(2)}$ in~$Q^{(1)}$, I have not written out a proof.}
$2$ in~$Q^{(1)}$.
It is the failure of this surjectivity that shows that the system
defined by~$Q$ is not formally integrable, which, of course,
implies that it is not involutive.
\end{remark}

\begin{remark}[Generality revisited]
\label{rem: Genrevisit}
In light of Corollary~\ref{cor: qidichotomy}, 
the study of (elliptic) quasi-integrable 
$\bbR^+{\cdot}\U(3)$-structures~$\bigl(J,[\eta]\bigr)$ 
can now be broken into two (overlapping) classes:  

The first class is defined by~$\lambda=0$ and is easy to understand.
In this case,~$J$ is integrable and it is easy now to argue that
these depend on~$14$ functions of~$6$ variables, locally.  To see
this, note that, because~$J$ is integrable, one can imagine constructing
these in two stages:  First, choose an integrable almost complex 
structure~$J$ on~$M$.  Since these are all locally equivalent
up to diffeomorphism, the general integrable almost complex structure
on~$M$ is easily seen to depend on~$6$ arbitrary functions of~$6$
variables (i.e., the generality of the diffeomorphism group in 
dimension~$6$).  Second, once~$J$ is fixed, choose a section~$[\eta]$ 
up to positive multiples of the bundle of positive~$(1,1)$-forms
for~$J$.  Since this projectivized bundle has real rank~$8$ over~$M$,
it follows that this choice depends on~$8$ functions of~$6$ variables.
Consequently, the general integrable conformal unitary structure
on~$M^6$ depends on~$14=6+8$ arbitrary functions of~$6$ variables.
(Of course, when one reduces modulo diffeomorphism, this comes
back down to~$8$ function of $6$ variables.)

The second class, defined by the first-order quasi-integrable conditions 
plus the~$6$ second-order conditions~\eqref{eq: QI2ndorderreln},
is more problematic.  If one lets~$Q''_1\subset Q^{(1)}$
denote the codimension~$6$ submanifold defined by these equations,
one does not know that~$Q''_1$ is formally integrable, let alone
involutive.  My attempts to check involutivity for~$Q''_1$ have, 
so far, been unsuccessful. Thus, it is hard to say how many functions 
of how many variables the general solution depends on. 

What one does know is that the tableau of~$Q''_1$ is isomorphic
to a proper subtableau of the system in~$J^2(\CU(M))$ that defines 
the conformally unitary structures with integrable underlying almost
complex structure and its last nonzero character is strictly
less than~$14$.  Consequently, generality of quasi-integrable
conformally unitary structures of the second type is strictly 
less than~$14$ functions of~$6$ variables.  

On the other hand, there is a lower bound on this generality 
that can be derived by a construction in the next subsubsection,
where a family depending on $10$ arbitrary functions of~$6$
variables is constructed by making use of some
linear algebra constructions involving $3$-forms on $6$-manifolds.
\end{remark}

\subsubsection{A construction}
I do not know an explicit construction 
of the general quasi-integrable confomally unitary structure.
However, there is a natural construction of a large class 
of such structures (depending, as will be seen, 
on $10$ functions of~$6$ variables) 
which I will now outline.  This construction, which
I found in 2000, has much in common with a construction 
of Hitchin's~\cite[\S5]{MR1871001}, and the reader 
may want to compare his treatment.

First, recall some linear algebra:%
\footnote{For a proof of these statements, see the Appendix.}  
For any $6$-dimensional
vector space~$V$ over~$\bbR$, the space~$A^3(V)=\Lambda^6(V^*)$ 
of alternating $3$-forms has dimension~$20$ and contains two
open~$\GL(V)$-orbits~$A^3_\pm(V)\subset A^3(V)$.  If~$v^i$
is a basis of~$V^*$, the open set~$A^3_+(V)$ is the $\GL(V)$-orbit
of
\be
\begin{aligned}
\phi_+ &= \Re\bigl((v^1+\iC v^2)\w(v^3+\iC v^4)\w(v^5+\iC v^6)\bigr)\\
       &= v^1\w v^3\w v^5-v^1\w v^4\w v^6-v^2\w v^3\w v^6-v^2\w v^4\w v^5
\end{aligned}
\ee
while the open set~$A^3_-(V)$ is the $\GL(V)$-orbit of
\be
\phi_- = v^1\w v^2\w v^3 + v^4\w v^5\w v^6. 
\ee
My interest will be in~$\phi_+$.

Let~$J_+:V\to V$ be the complex structure on~$V$ 
for which the complex-valued $1$-forms~$(v^1+\iC v^2)$, 
$(v^3+\iC v^4)$, and $(v^5+\iC v^6)$ are $J$-linear
and let~$\GL(V,J_+)\subset\GL(V)$ be the commuting subgroup of~$J_+$.
Of course, $\GL(V,J_+)$ is isomorphic to~$\GL(3,\bbC)$ and hence contains
a simple, normal subgroup~$\SL(V,J_+)\subset\GL(V,J_+)$ isomorphic
to~$\SL(3,\bbC)$ that consists of the elements of~$\GL(V,J_+)$ that
are complex unimodular.
 
If~$C\in\GL(V)$ anticommutes with~$J_+$ and satisfies~$\det(C)=-1$, 
then the $\GL(V)$-stabilizer of~$\phi_+$ is easily seen 
to be the $2$-component group~$\SL(V,J_+)\cup C{\cdot}\SL(V,J_+)$.  
In fact, $\SL(V,J_+)$ is the intersection of the stabilizer of~$\phi_+$
with~$\GL^+(V)$, the group of orientation-preserving linear transformations 
of~$V$. 

It follows that if~$V$ is endowed with an orientation and~$\phi$
is any element of~$A^3_+(V)$, then there is a unique 
element~$\ast\phi\in A^3_+(V)$ and complex structure~$J_\phi:V\to V$
such that~$\phi\w{\ast\phi}>0$ and such that~$\phi+\iC\,{\ast\phi}$ 
is a $(3,0)$-form for~$J_\phi$.  
Note the identity~$J_\phi^*(\phi)=\ast\phi$.
Note also that reversing the chosen orientation of~$V$ 
will replace~$(\ast\phi,J_\phi)$ by~$(-{\ast\phi},-J_\phi)$.
  
It is important to bear in mind
that the map~$\ast: A^3_+(V)\to A^3_+(V)$, 
though it is obviously smooth algebraic, 
is homogeneous of degree~$1$, and satisfies~$\ast{\ast\phi} = -\phi$, 
is \emph{not} the restriction to~$A^3_+(V)$ 
of a linear endomorphism of~$A^3(V)$. 

Of course, $J_\phi$ induces a decomposition of the complex-valued
alternating forms on~$V$ into types.  I will denote the associated
subspaces by~$A^{p,q}_\phi(V)$, as usual.  For my purposes, the
most important thing to note is that the space~$A^{p,p}_\phi(V)$
is the complexification of the real subspace~$H^p_\phi(V)$ consisting
of the real-valued forms in~$A^{p,p}_\phi(V)$.  Moreover, $H^p_\phi(V)$
contains an open convex cone~$H^p_\phi(V)^+$ consisting of the forms 
that are positive on all $J_\phi$-complex subspaces of dimension~$p$.

In particular, note that the normalized 
squaring map~$s:H^2_\phi(V)^+\to H^4_\phi(V)^+$
defined by~$s(\omega) = \frac12\omega^2$
is a diffeomorphism and hence has a smooth inverse
\be
\sigma:H^4_\phi(V)^+ \to H^2_\phi(V)^+.
\ee 
There is a hyperbolic analog of this square root:  Let~$H^2_\phi(V)^{r,s}$
denote the open set of nondegenerate $J_\phi$-Hermitian forms 
of type~$(r,s)$ (so that~$H^2_\phi(V)^+ = H^2_\phi(V)^{3,0}$).  Then
the normalized squaring map~$s:H^2_\phi(V)^{1,2}\to H^4_\phi(V)$ is
a diffeomorphism onto its (open) image, denoted~$H^4_\phi(V)^{1,2}$.
Accordingly, I extend the domain of~$\sigma$ to include~$H^4_\phi(V)^{1,2}$
and so that it inverts~$s$ on~$H^4_\phi(V)^{1,2}$.

Now let~$M^6$ be an oriented $6$-manifold and let~$\cA^{3}_+(M)$
denote the set of $3$-forms on~$M$ that, at each point are 
linearly equivalent to~$\phi_+$.  These forms are the sections 
of an open subset of~$A^3(TM)$ and hence are \emph{stable} in
Hitchin's sense~\cite{MR1871001}.  

For each~$\varphi\in\cA^3_+(M)$, 
there is associated a smooth~$3$-form~$\ast\varphi\in\cA^3_+(M)$ 
and an associated almost complex structure~$J_\varphi:TM\to TM$.  

Using the type decomposition with respect to~$J_\varphi$
(and the fact that~$\d^{-2,3}=0$), 
there exist a $(1,0)$-form~$\beta$ and a $(2,2)$-form~$\pi$ 
on~$M$ such that
\be
\label{eq: d30form}
\d(\varphi + \iC\,{\ast\varphi})
 = \overline{\beta} \w (\varphi + \iC\,{\ast\varphi})+\pi.
\ee
Of course~$\pbar(\varphi + \iC\,{\ast\varphi}) = 
\overline{\beta} \w (\varphi + \iC\,{\ast\varphi})$ 
and~$\d^{-1,2}(\varphi + \iC\,{\ast\varphi}) = \pi$.

Now let~$\varphi\in\cA^3_+(M)$ have the property that~$\ast\varphi$
is closed.  Then the left and right hand sides of~\eqref{eq: d30form}
are both real and, by the type decomposition, it follows that~$\beta=0$
and~$\pi = \bar\pi$ lies in~$\cH^4(M)$, the space of sections
of the bundle~$H^4_\varphi(TM)$.  

I will say that~$\varphi$ is \emph{positive definite} 
if~$\pi$ is a section of~$H^4_\varphi(TM)^+$ 
and that~$\varphi$ is \emph{positive indefinite} if~$\pi$
is a section of~$H^4_\varphi(TM)^{1,2}$.

Obviously, positivity (of either type) is an open condition 
on the $1$-jet of~$\varphi$.

If~$\varphi$ satisfies the condition that~$\d(\ast\varphi)=0$
and is positive (definite or indefinite), then there exists 
a unique~$(1,1)$-form~$\eta_\varphi$ 
satisfying~$\d\varphi = \pi = 2\,{\eta_\varphi}^2$
that lies in~$\cH^2(M)^+$ if~$\varphi$ is positive definite
and in~$\cH^2(M)^{1,2}$ if $\varphi$ is positive indefinite.

\begin{proposition}[Quasi-integrability]
\label{prop: qiconfclosed}
Let~$\varphi\in\cA^3_+(M)$ satisfy~$\d(\ast\varphi)=0$. 
If~$\varphi$ is positive definite, 
then the $\U(3)$-structure~$(J_\varphi,\eta_\varphi)$
is strictly quasi-integrable.  If~$\varphi$ is positive indefinite,
then the~$\U(1,2)$-structure~$(J_\varphi,\eta_\varphi)$
is strictly quasi-integrable.
\end{proposition}

\begin{proof}
First, assume that~$\varphi$ is positive definite.
The claimed result is local, so let~$\alpha = (\alpha^i):TU\to\bbC^3$ 
be a local $(J_\varphi,\eta_\varphi)$-unitary coframing.  Then
\be
\eta_\varphi = {\ts\frac\iC2}\,{}^t\alpha\w\bar\alpha,
\ee
so that
\be
\label{eq: piasalpha}
\pi = 2\,{\eta_\varphi}^2 = 
\bpm\alpha^2\w\alpha^3&\alpha^3\w\alpha^1&\alpha^1\w\alpha^2\epm\w
\bpm\overline{\alpha^2\w\alpha^3}\strut\\
    \overline{\alpha^3\w\alpha^1}\strut\\
    \overline{\alpha^1\w\alpha^2}\strut
\epm.
\ee
Moreover, there exists a nonvanishing complex function~$F$ on~$U$ such that
\be
\label{eq: newFdefined}
\varphi + \iC\,{\ast\varphi} = F^{-1}\,\alpha^1\w\alpha^2\w\alpha^3,
\ee
implying that
\be
\pi = \d^{-1,2}(\varphi + \iC\,{\ast\varphi})
= F^{-1}\,\bpm\alpha^2\w\alpha^3&\alpha^3\w\alpha^1&\alpha^1\w\alpha^2\epm
N(\alpha)
\bpm\overline{\alpha^2\w\alpha^3}\strut\\
    \overline{\alpha^3\w\alpha^1}\strut\\
    \overline{\alpha^1\w\alpha^2}\strut
\epm.
\ee
Comparison with~\eqref{eq: piasalpha} now yields that~$N(\alpha) = F\,I_3$.  
In particular, $N(\alpha)$ is a nonzero multiple of the identity, 
which is what needed to be shown.

The proof in the positive indefinite case is completely similar.  One
merely has to note that the formulae become
\be
\eta_\varphi = {\ts\frac\iC2}\,
\bigl(\alpha^1\w\overline{\alpha^1} - \alpha^2\w\overline{\alpha^2}
     - \alpha^3\w\overline{\alpha^3} \bigr)
\ee
so that
\be
\pi = 2\,{\eta_\varphi}^2 = 
\bpm\alpha^2\w\alpha^3&-\alpha^3\w\alpha^1&-\alpha^1\w\alpha^2\epm\w
\bpm\overline{\alpha^2\w\alpha^3}\strut\\
    \overline{\alpha^3\w\alpha^1}\strut\\
    \overline{\alpha^1\w\alpha^2}\strut
\epm.
\ee
The rest of the proof is entirely similar to the positive definite case.
\end{proof}

\begin{remark}[Further information and a characterization]
The calculations made in the course of the proof can be carried
a little further:  Assuming that~$U$ is simply connected, one can,
by multiplying~$\alpha$ by a unimodular complex function, arrange
that~$F$ is real and positive.  

For example, imposing this condition in the positive definite case 
specifies~$\alpha$ up to replacement of the form~$\alpha\mapsto g\alpha$ 
where~$g:U\to\SU(3)$ is smooth.  Since~$N(\alpha) = F\,I_3$, 
it follows from the definitions of~$\omega$ and of~$\psi$ 
that~$\omega = F^2\,\eta_\varphi$ and
that~$\psi = F^3\,\alpha^1\w\alpha^2\w\alpha^3 
= F^4\,(\varphi + \iC\,{\ast\varphi})$.

Using this information, it is easy to see that the strictly
quasi-integrable structures~$\bigl(J,\eta\bigr)$ 
constructed in Proposition~\ref{prop: qiconfclosed} 
are characterized by the condition
that the canonical conformal curvature form~$\d\rho$ vanishes identically,
or, equivalently, in terms of the almost complex structure~$J$,
that there be a positive function~$F$ on~$M$
such that~$F^{-4}\,\omega(J)^2$ is closed.
(Such an~$F$, if it exists, is clearly unique up to constant multiples.)
\end{remark}

\begin{remark}[Generality]
In Cartan's sense of generality, the closed $3$-forms in dimension~$6$
depend on~$10$ functions of~$6$ variables.  It is an open condition
(of order~$0$) on a $3$-form~$\varphi$ that it lie in~$\cA^3_+(M)$
and it is a further open condition (of order~$1$) on~$\varphi$ that
$\d\varphi$ lie in~$\cH^4(M)^+\cup\cH^4(M)^{1,2}$.  Thus, one can
say that the local positive (definite or indefinite) 
$3$-forms~$\phi\in \cA^3_+(M)$ satisfying~$\d(\ast\phi)=0$ 
depend on $10$ functions of~$6$ variables.

These conditions are obviously invariant under diffeomorphisms
in dimension~$6$, so it makes sense to say that germs of solutions
modulo diffeomorphism depend on~$10-6=4$ arbitrary functions 
of~$6$ variables.  (Strictly speaking, in order to make this count
work, one needs to observe that the group of symmetries of any
particular~$\phi$ satisfying these conditions is the group of 
symmetries of the associated~$\SU(3)$- or~$\SU(1,2)$-structure
and hence is finite dimensional.)  This count can be made more
rigorous (and verified) by appealing to Cartan's theory of generality, 
but I will not do this analysis here.
\end{remark}

\appendix

\section{$3$-forms in Dimension~$6$}

In this appendix, I will supply a proof of the following
normal form result, whose analog over the complex field is well-known 
and due to Reichel~\cite{Reichel1907}, but whose complete proof 
over the real field does not appear to be easy to locate 
in the literature.%
\footnote{Probably, this is due to my inability to read Russian.
If any reader can point me to a proof in the literature, 
I'll be grateful.}  
For applications in this article, 
the important case is the second normal form.

\begin{proposition}[Normal forms] 
\label{prop: normalforms}
  Let $V$ be a real vector space of dimension~$6$ over~$\bbR$ and let
$\phi\in\Lambda^3(V^*)$ be any element.  Then there exists a basis~
$e^1,\ldots,e^6$ of~$V^*$ so that~$\phi$ is equal to one of the following 
\begin{enumerate}
\item $e^1\w e^2\w e^3 + e^4\w e^5\w e^6$,
\item $ e^1\w e^3\w e^5 - e^1\w e^4\w e^6
       -e^2\w e^3\w e^6 - e^2\w e^4\w e^5$,
\item $e^1\w e^5\w e^6 + e^2\w e^6\w e^4
       +e^3\w e^4\w e^5$,
\item $e^1\w e^2\w e^5 + e^3\w e^4\w e^5$,
\item $e^1\w e^2\w e^3$, or
\item $0$.
\end{enumerate}
Moreover, these six forms are mutually inequivalent.  The~$\GL(V)$-orbits
of the first two $3$-forms are open in~$\Lambda^3(V^*)$, the $\GL(V)$-orbit
of the third $3$-form is a hypersurface in~$\Lambda^3(V^*)$, and the
$\GL(V)$-orbits of the remaining forms are of higher codimension.
\end{proposition}

\begin{proof}
The proof will be a series of steps, beginning with a sort of zeroth step
to take care of certain special cases.  A form~$\phi\in\Lambda^3(V^*)$ 
will be said to be {\it degenerate} if there is a non-zero vector~$v\in V$
that satisfies~$v\lhk\phi=0$.  Equivalently, $\phi$ is degenerate
if there exists a subspace~$v^\perp\subset V^*$ of rank~$5$ so that~$\phi$ 
lies in~$\Lambda^3(v^\perp)$.   In such a case, using the canonical 
isomorphism
\be
\Lambda^3(v^\perp)=\Lambda^2\bigl((v^\perp)^*\bigr)\otimes\Lambda^5(v^\perp)
\ee
together with the well-known classification of 2-forms, it follows
that either $\phi=0$ or else is there is a basis~$e_1,\ldots,e_5$ of
$(v^\perp)^*$, with dual basis~$e^1,\ldots,e^5$ of~$v^\perp$, so that
$\phi$ is either
\be
(e_1\w e_2 + e_3\w e_4)\otimes e^1\w e^2\w e^3\w e^4\w e^5
 = e^1\w e^2\w e^5 + e^3\w e^4\w e^5 
\ee
or
\be
  (e_4\w e_5)\otimes e^1\w e^2\w e^3\w e^4\w e^5 
  = e^1\w e^2\w e^3.
\ee
Thus, the degenerate cases are accounted for by the last three types
listed in the Proposition.

It is useful to note that $\phi$ is degenerate if and only if~$\phi$ 
admits a linear factor, i.e., is of the form~$\phi=\alpha\w\beta$ for 
some~$\alpha\in V^*$ and $\beta\in\Lambda^2(V^*)$. To see this, note that
each of the three degenerate types has at least one linear factor. Conversely, 
given the linear factor~$\alpha$, the $2$-form~$\beta$ can actually be 
regarded as a well-defined element of~$\Lambda^2\bigl(V^*/(\bbR\,\alpha)\bigr)$.
Since the quotient space~$V^*/(\bbR\,\alpha)$ has dimension~$5$, it follows
that~$\beta$ is either zero, decomposable, or the sum of two decomposable  
terms, again yielding forms of the last three types listed in the
Proposition.  

Henceforth, it will be assumed that~$\phi$ is nondegenerate, i.e., that
$v\lhk\phi\not=0$ for all nonzero~$v\in V$ and also that~$\phi$ has no
linear divisors, i.e., that~$\alpha\w\phi\not=0$ for all nonzero~$\alpha\in V$.
It remains to show that a nondegenerate element of~$\Lambda^3(V^*)$ can be
put into one of the first three forms listed in the Proposition.

It is convenient to start with what is apparently a special case.  Let~
$e_1,\ldots,e_6$ be a basis of~$V$, with dual basis~$e^1,\dots,e^6$ of~$V^*$.
Set
\be
\phi_0 = e^1\w e^2\w e^3 + e^4\w e^5\w e^6.
\ee
I claim that the $\GL(V)$-orbit of~$\phi_0$ is open in~$\Lambda^3(V^*)$
and will establish this by showing that the dimension of the stabilizer
\be
G_0 = \left\{\ g\in\GL(V)\ \vrule\ g^*(\phi_0) = \phi_0\ \right\}
\ee
is~16, so that the orbit~$\GL(V){\cdot}\phi_0 = \GL(V)/G_0$ has dimension
$36-16 = 20$, which is the dimension of~$\Lambda^3(V^*)$.

Let~$C_0 = \{\ v\in V\ \vrule\ (v\lhk \phi)^2 = 0\ \}$.  Computation shows
that~$C_0 = P_0^+\cup P_0^-$ where~$P_0^+ = \text{span}\{e_1,e_2,e_3\}$ and
$P_0^- = \text{span}\{e_4,e_5,e_6\}$.  Since the elements of~$G_0$ must 
preserve~$C_0$, it follows that they permute the subspaces~$P_0^\pm$.  
Let~$G'_0\subset G_0$ be the subgroup consisting of elements that preserve the 
subspaces~$P_0^+$ and~$P_0^-$.  Since the linear transformation~$c:V\to V$
defined by
\be
c(e_i) = \begin{cases} e_{i+3}& i=1,2,3, \\ 
                       e_{i-3}& i=4,5,6. \\ 
         \end{cases}
\ee
lies in~$G_0$ but not in~$G'_0$, it follows that~$G_0=G'_0\cup G'_0{\cdot}c$.
Since~$\phi_0$ pulls back to each of~$P_0^\pm$ to be a volume form,
it follows that~$G'_0 = \SL(P_0^+)\times\SL(P_0^-)\simeq\SL(3,\bbR)\times
\SL(3,\bbR)$.  Thus,~$\dim G_0 = 16$, as desired.

Now, for any~$\phi\in\Lambda^3(V)$ (degenerate or not), define the map
\be
J_\phi : V\to \Lambda^5(V^*) = V\otimes \Lambda^6(V^*)
\ee
by~$J_\phi(v) = (v\lhk\phi)\w\phi$, and regard~$J_\phi$ as
an element of~$\End(V)\otimes \Lambda^6(V^*)$.  Then the 
map~$\phi\mapsto J_\phi$ is a quadratic polynomial map from~$\Lambda^3(V^*)$
to~$\End(V)\otimes \Lambda^6(V^*)$ that is equivariant with respect to the
natural actions of~$\GL(V)$ on the domain and range.

 For simplicity of notation, write~$J_0$ for~$J_{\phi_0}$.  
Inspection shows that
\be
\begin{aligned}
J_0 &= \left( e_1\otimes e^1 + e_2\otimes e^2 + e_3\otimes e^3
            -e_4\otimes e^4 - e_5\otimes e^5 - e_6\otimes e^6 \right)\\
&\qquad\qquad\qquad\qquad\qquad \otimes e^1\w e^2\w e^3\w e^4\w e^5\w e^6.\\
\end{aligned}
\ee
In particular, it follows that ~$\tr J_0 = 0\in\Lambda^6(V^*)$.  Because of
the $\GL(V)$-equivariance already mentioned, it follows that~$\tr J_\phi=0$
for all~$\phi$ in the orbit $\GL(V){\cdot}\phi_0$.  Since this latter orbit is
open and since~$J$ is a polynomial mapping, it follows that~$\tr J_\phi=0$
for all~$\phi\in\Lambda^3(V^*)$.  

Now consider~${J_\phi}^2\in\End(V)\otimes S^2\bigl(\Lambda^6(V^*)\bigr)$.
By the above formula, 
\be
\begin{aligned}
{J_0}^2 &= \left( e_1\otimes e^1 + e_2\otimes e^2 + e_3\otimes e^3
            +e_4\otimes e^4 + e_5\otimes e^5 + e_6\otimes e^6 \right)\\
&\qquad\qquad\qquad\qquad\qquad \otimes 
\left(e^1\w e^2\w e^3\w e^4\w e^5\w e^6\right)^2.\\
\end{aligned}
\ee
In particular, 
$\tr {J_0}^2 = 6\,\left(e^1\w e^2\w e^3\w e^4\w e^5\w e^6\right)^2$, so that
\be
{J_0}^2  - {\ts\frac16}\,\id_V\otimes \tr \left({J_0}^2\right) = 0.
\ee
This implies that the $\GL(V)$-equivariant quartic polynomial map
from~$\Lambda^3(V^*)$ to~$\End(V)\otimes S^2\bigl(\Lambda^6(V^*)\bigr)$ 
defined by
\be
\phi\mapsto {J_\phi}^2 -{\ts\frac16}\,\id_V\otimes\tr\left({J_\phi}^2\right)
\ee
vanishes on the~$\GL(V)$-orbit of~$\phi_0$, which is open.  Thus, it
follows that the identity
\be
{J_\phi}^2  = {\ts\frac16}\,\id_V\otimes\tr\left({J_\phi}^2\right)
\ee
holds for all~$\phi\in\Lambda^3(V^*)$.

The remainder of the proof will divide into cases according to the
$\GL(V)$-orbit of~$\tr\left({J_\phi}^2\right)$ 
in the 1-dimensional vector space~$S^2\bigl(\Lambda^6(V^*)\bigr)$.  
There are three such orbits and they can naturally be designated as 
`positive', `negative', and `zero'. 

Suppose first that~$\tr\left({J_\phi}^2\right)$ be `positive', i.e., that it
can be written in the form~$6\,\Omega^2$ for some~$\Omega\in\Lambda^6(V^*)$,
which is unique up to sign.  By the above identity, it then follows that
${J_\phi}^2(v) = v\otimes \Omega^2$ for all~$v\in V$. Define two 
`eigenspaces' of~$J_\phi$ by
\be
P_\phi^\pm = \left\{\ v\in V\ \vrule\ J_\phi(v) = \pm v\lhk\Omega
= \pm v\otimes \Omega\ \right\}.
\ee
Because of the identity~${J_\phi}^2(v) = v\otimes \Omega^2$, it follows
that~$V = P_\phi^-\oplus P_\phi^+$. Since~$\tr J_\phi=0$, it follows
that~$\dim P_\phi^- = \dim P_\phi^+ = 3$.  Note that replacing~$\Omega$
by~$-\Omega$ will simply switch these two spaces.

Now, this splitting of~$V$ implies a splitting of~$\Lambda^3(V^*)$ as
\be
\begin{aligned}
\Lambda^3(V^*) &= \Lambda^3\bigl((P_\phi^+)^*\bigr)
\oplus \Lambda^2\bigl((P_\phi^+)^*\bigr)\w (P_\phi^-)^*\\
&\qquad{}\oplus (P_\phi^+)^*\w \Lambda^2\bigl((P_\phi^-)^*\bigr)
         \oplus\Lambda^3\bigl((P_\phi^-)^*\bigr)
\end{aligned}
\ee
and hence a corresponding splitting of~$\phi$ into four terms as
\be
\phi = \phi^{+++} + \phi^{++-} + \phi^{+--} + \phi^{---}.
\ee

I will now argue that~$\phi^{++-}=\phi^{+--}=0$, which will show that
$\phi$ has the form~(1), since neither $\phi^{+++}$ nor $\phi^{---}$
can vanish for such a~$\phi$ if it be nondegenerate.  
To see this, let~$v^+\in P_\phi^+$
and~$v^-\in P_\phi^-$ be arbitrary.  By definition, $(v^+\lhk\phi)\w\phi
= v^+\lhk\Omega$ while $(v^-\lhk\phi)\w\phi= -v^-\lhk\Omega$.  Taking
left hooks then yields the identities
\be
\begin{aligned}
 v^-\lhk(v^+\lhk\Omega) 
&= \bigl(v^-\lhk(v^+\lhk\phi)\bigr)\w\phi + (v^+\lhk\phi)\w(v^-\lhk\phi)\\
-v^+\lhk(v^-\lhk\Omega) 
&= \bigl(v^+\lhk(v^-\lhk\phi)\bigr)\w\phi + (v^-\lhk\phi)\w(v^+\lhk\phi)\\
\end{aligned}
\ee
Now, subtracting these identities 
yields~$2\bigl(v^-\lhk(v^+\lhk\phi)\bigr)\w\phi=0$.  If there did exist~$v^+$
and~$v^-$ so that~$v^-\lhk(v^+\lhk\phi)$ were nonzero, then this equation 
would imply that~$\phi$ had a linear divisor, contradicting the nondegeneracy 
assumption. Thus, it
must be true that $v^-\lhk(v^+\lhk\phi)=0$ for all $v^+\in P_\phi^+$
and~$v^-\in P_\phi^-$.  However, this is equivalent to the condition that~
$\phi^{++-}=\phi^{+--}=0$, as desired.  It follows that there is a basis~
$e^1,\ldots, e^6$ of~$V^*$ so that~$\phi^{+++}=e^1\w e^2\w e^3$ while~
$\phi^{---}=e^4\w e^5\w e^6$.  Thus, $\phi$ is of the form~$(1)$.

Next, suppose that~$\tr\left({J_\phi}^2\right)$ be `negative', i.e., that it
can be written in the form~$-6\,\Omega^2$ for some~$\Omega\in\Lambda^6(V^*)$,
which is unique up to sign.  By the above identity, it then follows that
${J_\phi}^2(v) = -v\otimes \Omega^2$ for all~$v\in V$. Define two `complex
eigenspaces' of~$J_\phi$ by
\be
P_\phi^\pm = \left\{\ v\in V^\bbC \ \vrule\ J_\phi(v) = \pm i\, v\lhk\Omega
= \pm i\, v\otimes \Omega\ \right\}.
\ee
Note that since~$\Omega$ is real, these two subspaces of~$V^\bbC=V\otimes\bbC$
are conjugate.
Reasoning as in the positive case now yields that there must
exist a basis~$e^1,\ldots, e^6$ of~$V^*$ so that~$\phi=\phi^{+++}+\phi^{---}$
where~$\phi^{+++} = {\frac12}\,(e^1+\iC\,e^2)\w(e^3+\iC\,e^4)\w(e^5+\iC\,e^6)$ 
while~$\phi^{---} = \overline{\phi^{+++}} =
{\frac12}\,(e^1-\iC\,e^2)\w(e^3-\iC\,e^4)\w(e^5-\iC\,e^6)$.  
This yields the form~$(2)$.

Finally, suppose that $\tr\left({J_\phi}^2\right)=0$, so that, by
the identity above, ${J_\phi}^2 = 0$.  Let~$K_\phi\subset V$ be the kernel
of~$J_\phi$ and let~$I_\phi\subset V$ be such that the image of~$J_\phi$
is $I_\phi\otimes\Lambda^6(V^*)\subset V\otimes\Lambda^6(V^*)=\Lambda^5(V^*)$.
Then $\dim K_\phi + \dim I_\phi = 6$ and $I_\phi\subset K_\phi$ 
(since~${J_\phi}^2=0$), so $\dim K_\phi\ge 3$. 

Now, for $v\in K_\phi$, the identity~$(v\lhk\phi)\w\phi=0$ holds.  It follows
that, for $v_1,v_2\in K_\phi$, the identity
\be
\bigl(v_1\lhk(v_2\lhk\phi)\bigr)\w\phi + (v_2\lhk\phi)\w(v_1\lhk\phi)=0
\ee
must also hold.  However, the first term in this expression is skewsymmetric
in~$v_1,v_2$ while the second term is symmetric in~$v_1,v_2$.  It follows
that each term must vanish separately, i.e., that
\be
\bigl(v_1\lhk(v_2\lhk\phi)\bigr)\w\phi = (v_1\lhk\phi)\w(v_2\lhk\phi)=0
\ee
for all $v_1,v_2\in K_\phi$.  Since~$\phi$ is nondegenerate, it has no
linear divisors, so the first of these equations implies that 
\be
v_1\lhk(v_2\lhk\phi) = 0
\ee
for all $v_1,v_2\in K_\phi$.  This implies that~$\phi$ must lie
in the subspace~$\Lambda^2({K_\phi}^\perp)\w V^*\subset\Lambda^3(V^*)$.  Now
if~$\dim K_\phi$ were greater than~$4$, then~$\Lambda^2({K_\phi}^\perp)$
would vanish, which is absurd.  If~$\dim K_\phi$ were equal to~$4$,
then~$\Lambda^2({K_\phi}^\perp)$ would be spanned by a single decomposable
$2$-form, which would imply that~$\phi$ has two linearly independent
divisors, which is also impossible for a nondegenerate form.  Thus~$
\dim K_\phi=3$ is the only possibility.  Let~$e^4,e^5,e^6$ be a basis
of~${K_\phi}^\perp$ and write
\be
\phi = e^1\w e^5\w e^6 + e^2\w e^6\w e^4 + e^3\w e^4\w e^5 
         + c\,e^4\w e^5\w e^6
\ee
for some~$e^1,e^2,e^3\in V^*$ and some constant~$c\in\bbR$.  
The $1$-forms~$e^1,e^2,e^3,e^4,e^5,e^6$ must be linearly independent 
and hence a basis of~$V^*$ or else~$\phi$ would be degenerate. 
Now replacing~$e^1$ by $e^1 + c\,e^4$ in the above basis yields a basis
in which
\be
\phi = e^1\w e^5\w e^6 + e^2\w e^6\w e^4 + e^3\w e^4\w e^5 \,,
\ee 
showing that~$\phi$ lies in the orbit of the form~(3), as was desired.
\end{proof}

\begin{remark}[The stabilizer groups]
The proof of Proposition~\ref{prop: normalforms} 
indicates the stablizer subgroup~$G_\phi\subset\GL(V)$ 
for each of the normal forms.  

If~$\phi$ is of type~$(1)$, 
then~$G_\phi = G_\phi'\cup G_\phi'{\cdot}c$, where~$G_\phi'$
is the index~$2$ subgroup of~$G_\phi$ 
that preserves the two $3$-planes~$P_\phi^\pm$
and their volume forms while~$c:V\to V$ is a linear map that exchanges
these two planes and their volume forms.  Thus, there is an exact
sequence
\be
0\longrightarrow \SL(P_\phi^+)\times \SL(P_\phi^-)
 \longrightarrow G_\phi
 \longrightarrow \bbZ_2
 \longrightarrow 0.
\ee

If~$\phi$ is of type~$(2)$, then $G_\phi = G_\phi'\cup G_\phi'{\cdot}c$
where~$G_\phi'$ is the index two subgroup that preserves a complex structure
and complex volume form~$\omega = 2\phi^{+++}$ on~$V$ (i.e., 
$\omega$ is a decomposable element of~$\Lambda^3(V\otimes \bbC)$), 
while~$c:V\to V$ is a conjugate linear mapping of~$V$ to itself that 
takes~$\omega$ to ~$\overline{\omega}$. In fact, $\omega$ can be chosen 
so that~$\phi = {\frac12}(\omega + \overline{\omega})$.
Thus, there is an exact sequence
\be
0\longrightarrow \SL(V^\bbC)
 \longrightarrow G_\phi
 \longrightarrow \bbZ_2
 \longrightarrow 0.
\ee
The important thing to note is that~$G_\phi\cap \GL^+(V) = G_\phi'
\simeq \SL(V^\bbC)$.  Thus, the orientation-preserving stabilizer 
of~$\phi$ preserves a unique complex structure for which~$\phi$ is
the real part of a $(3,0)$-form.

The other stabilizers will not be important here, 
so we leave those to the reader.
\end{remark}

\begin{remark}[The symplectic version]
Various authors~\cite{MR2002657,MR0725059} 
have considered and solved a symplectic version 
of the above normal form problem:  

Namely, consider a $6$-dimensional
real vector space~$V$ endowed with a symplectic (i.e., nondegenerate)
$2$-form~$\omega\in\Lambda^2(V^*)$  and let~$\Symp(\omega)\subset\GL(V)$
denote the stabilizer subgroup of~$\omega$.  

As a representation of~$\Symp(\omega)$, 
the vector space~$\Lambda^3(V^*)$ is reducible,
being expressible as a direct sum~$\Lambda^3(V^*) = \omega{\w}V^*
\oplus \Lambda^3_0(V^*)$, where~$\Lambda^3_0(V^*)$ 
is the space of $3$-forms~$\phi$ on~$V$ that satisfy~$\omega\w\phi=0$.
These two summands are irreducible $\Symp(\omega)$ modules and it
is an interesting problem to classify the orbits of~$\Symp(\omega)$
acting on~$\Lambda^3_0(V^*)$.  

Proposition~\ref{prop: normalforms} 
provides a simple solution to this classification problem,
originally completed by Banos~\cite{MR2002657}.

For example, suppose that~$\phi\in\Lambda^3_0(V^*)$ has the normal
form~(1) and let~$e^1,\ldots,e^6$ be a basis of~$V^*$ such that
\be
\phi = e^1\w e^2\w e^3 + e^4\w e^5\w e^6.
\ee
The condition~$\omega\w\phi=0$ is equivalent to
the condition that there be~$a_{ij}$ such that
\be
\omega = a_{ij}\,e^i\w e^{j+3}
\ee
(where~$i$ and~$j$ run from~$1$ to~$3$ in the summation).  
Since~$\omega$ is non-degenerate,~$\det(a)$ is nonvanishing.
Thus, writing~$a = \lambda\,b$ where~$\det(b)=1$ and~$\lambda\not=0$, 
one can make a unimodular basis change in~$e^1,e^2,e^3$ so that
\be
\begin{aligned}
\omega &= \lambda\bigl(e^1\w e^4 + e^2\w e^5 + e^3\w e^6\bigr)\\
\phi &= e^1\w e^2\w e^3 + e^4\w e^5\w e^6\,.
\end{aligned}
\ee
Exchanging~$e^1,e^2,e^3$ with~$e^4,e^5,e^6$ if necessary, 
it can be assumed that~$\lambda$ is positive 
and then an overall scale change 
puts the pair $(\omega,\phi)$ into the form
\be
\label{eq: splitform}
\begin{aligned}
\omega &= e^1\w e^4 + e^2\w e^5 + e^3\w e^6\\
\phi &= \mu\bigl(e^1\w e^2\w e^3 + e^4\w e^5\w e^6\bigr)
\end{aligned}
\ee
with~$\mu>0$.  (Obviously, one cannot get rid of the factor~$\mu$ entirely.)

Similarly, if~$\phi$ is of type~$(2)$, 
then there are~$\zeta^1,\zeta^2,\zeta^3\in \bbC\otimes V^*$ such that
\be
\phi = {\ts\frac12}\bigl(\zeta^1\w\zeta^2\w\zeta^3 
             + \overline{\zeta^1\w\zeta^2\w\zeta^3}\,\bigr)\,.
\ee
The condition that~$\omega\w\phi=0$ is then equivalent to the condition
that
\be
\omega = {\ts\frac\iC2} a_{i\bj}\,\zeta^j\w\overline{\zeta^i}
\ee
for some Hermitian symmetric matrix~$a = (a_{i\bj}) = {}^t\bar a$
with nonvanishing determinant.  
Exchanging~$\zeta^i$ for~$\overline{\zeta^i}$ if necessary, 
it can be supposed that~$\det(a)=\lambda^3>0$ for some~$\lambda>0$.
By making a complex unimodular basis change in the~$\zeta^i$,
it can be further supposed that $a = \lambda\diag(1,\pm1,\pm1)$, 
reducing (after scaling) to the two possible normal forms
\be
\begin{aligned}
\omega &= {\ts\frac\iC2}\bigl(\zeta^1\w\overline{\zeta^1}
          \pm\zeta^2\w\overline{\zeta^2}\pm\zeta^3\w\overline{\zeta^3}\bigr) \\
\phi &= \mu \Re(\zeta^1\w\zeta^2\w\zeta^3)
\end{aligned}
\ee
for some~$\mu>0$.  
Alternatively, writing~$\zeta^1 = e^1+\iC e^4$, $\zeta^2 = e^2\pm\iC e^5$,
and $\zeta^3 = e^3\pm\iC e^6$ for a real basis~$e^i$ of~$V^*$, 
one can write this in the normal forms
\be
\label{eq: unitaryforms}
\begin{aligned}
\omega &= e^1\w e^4 + e^2\w e^5 + e^3\w e^6 \\
\phi &= \mu\bigl(e^1\w e^2\w e^3  - e^1\w e^5\w e^6 
                   \mp e^2\w e^6\w e^4 \mp e^3\w e^4\w e^5\bigr).
\end{aligned}
\ee

If~$\phi$ is of type~$(3)$, then there is a basis~$e^i$ of~$V^*$
such that
\be
\phi =  e^1\w e^5\w e^6 + e^2\w e^6\w e^4 + e^3\w e^4\w e^5.
\ee
The condition~$\omega\w\phi=0$ then implies that
\be
\omega = a_{ij}\,e^i\w e^{j+3} + b_{ij}\,e^{i+3}\w e^{j+3}
\ee
where~$a = {}^ta = (a_{ij})$ has~$\det(a)\not=0$ and~$b = -{}^tb = (b_{ij})$
is arbitrary.  Again, by making a basis change in the~$e^i$ that
preserves the form of~$\phi$, one can suppose that~$a = \diag(1,\pm1,\pm1)$.
Then by replacing~$e^i$ by~$e^i+s^i_j\,e^{j+3}$ for appropriate~$s^i_j$
satisfying~$s^i_i=0$, one can get rid of the~$b_{ij}$, leaving the
two normal forms
\be
\begin{aligned}
\omega &=  e^1\w e^4 \pm e^2\w e^5 \pm e^3\w e^6\\
\phi &= e^1\w e^5\w e^6 + e^2\w e^6\w e^4 + e^3\w e^4\w e^5.
\end{aligned}
\ee
or, alternatively,
\be
\label{eq: orthogonalforms}
\begin{aligned}
\omega &=  e^1\w e^4 + e^2\w e^5 + e^3\w e^6\\
\phi &= e^1\w e^5\w e^6 \pm e^2\w e^6\w e^4 \pm e^3\w e^4\w e^5.
\end{aligned}
\ee

By entirely similar arguments, one sees that when~$\phi$ is
of type~$(4)$, one can find a basis~$e^i$ of~$V^*$ such that
\be
\label{eq: degenerateforms}
\begin{aligned}
\omega &=  e^1\w e^4 + e^2\w e^5 + e^3\w e^6\\
\phi &= (e^1\w e^2 \pm e^4\w e^5)\w e^3,
\end{aligned}
\ee
(the two signs give inequivalent normal forms);
when~$\phi$ is of type~$(5)$, one can arrange
\be
\label{eq: lagrangianform}
\begin{aligned}
\omega &=  e^1\w e^4 + e^2\w e^5 + e^3\w e^6\\
\phi &= e^1\w e^2\w e^3;
\end{aligned}
\ee
and when~$\phi$ is of type~$(6)$, one can arrange
\be
\label{eq: nullform}
\begin{aligned}
\omega &=  e^1\w e^4 + e^2\w e^5 + e^3\w e^6\\
\phi &= 0.
\end{aligned}
\ee
This completes the list of normal forms.  

This analysis provides the stabilizer subgroups~$G(\omega,\phi)$ 
of the various normal forms~$(\omega,\phi)$.  
For example, in the case of~\eqref{eq: splitform},
the stabilizer subgroup is isomorphic to~$\SL(3,\bbR)$,
in the case of~\eqref{eq: unitaryforms} 
the stabilizer subgroup is isomorphic to~$\SU(3)$
if the upper sign is taken and~$\SU(1,2)$ if the lower sign is taken,
while in the case of~\eqref{eq: orthogonalforms}, the stabilizer
subgroup is isomorphic to the semidirect product of~$\bbR^5$ with
either $\SO(3)$ (upper sign) or $\SO(1,2)$ (lower sign).

\end{remark}

\bibliographystyle{hamsplain}

\providecommand{\bysame}{\leavevmode\hbox to3em{\hrulefill}\thinspace}

\end{document}